\documentclass[11pt, reqno]{amsart}
\usepackage{epsfig}

\usepackage{amstext}
\usepackage{xspace}
\usepackage{amsfonts}
\usepackage{amsmath}
\usepackage{amssymb}
\usepackage{amstext}
\usepackage{amsthm}       
\usepackage{xspace}
\usepackage[active]{srcltx}

\newcommand{\CC}{\D{C}}
\newcommand{\EE}{\mathbb E}

\newcommand{\R}{\mathbb R}

\newcommand{\N}{\mathbb N}

\newcommand{\comp}{\mbox{\scriptsize  $\circ$}}
\newcommand{\eps}{\varepsilon}

\newcommand{\ucv}{\rightrightarrows}

\newcommand{\A}{\mathcal{A}}

\newcommand{\oo}{\omega}
\newcommand{\Oo}{\Omega}
\newcommand{\de}{\delta}
\newcommand{\p}{\partial}

\newcommand{\F}{\mathcal{F}}

\newcommand{\M}{\mathcal{M}}
\renewcommand{\S}{\mathcal{S}}

\newcommand{\Tp}{T^+_t}
\newcommand{\Tm}{T^-_t}

\newcommand{\PP}{\mathbb{P}}
\newcommand{\tagliato}{$\kern-5.5 mm -$}
\newcommand{\tagliat}{$\kern-6.5 mm -$}

\newcommand{\cchi}{\mbox{\large $\chi$}}
\newcommand{\abra}[1]{(\ref{#1})}

\newcommand{\D}[1]{\mbox{\rm #1}}
\newcommand{\dd}{\D{d}}

\newtheorem{teorema}{Theorem}[section]
\newtheorem{prop}[teorema]{Proposition}
\newtheorem{lemma}[teorema]{Lemma}
\newtheorem{definition}[teorema]{Definition}
\newtheorem{cor}[teorema]{Corollary}

\newtheorem{guess}[teorema]{Remark}
\newtheorem{example}[teorema]{Example}

\newenvironment{dimo}{{\bf\noindent Proof.}}{\qed}
\newenvironment{oss}{\begin{guess} \begin{rm}}{\end{rm} \end{guess}}

\begin{document}

\title{Existence and regularity\\  of strict critical subsolutions\\
in the stationary ergodic setting}
\author{Andrea Davini \and Antonio Siconolfi}
\address{Dip. di Matematica, Universit\`a di Roma ``La Sapienza",
P.le Aldo Moro 2, 00185 Roma, Italy}
\email{davini@mat.uniroma1.it,\quad siconolf@mat.uniroma1.it }

\begin{abstract}
We prove that any continuous and convex stationary ergodic Hamiltonian 
admits  critical  subsolutions, which are strict outside
the random Aubry set. They   make up, in addition, a dense subset
of all critical subsolutions  with respect to a suitable metric.
If the Hamiltonian is additionally assumed of Tonelli type,  then
there exist strict subsolutions   of class $\CC^{1,1}$ in $\R^N$.
The proofs are based on the use of Lax--Oleinik semigroups and
their regularizing properties in the stationary ergodic
environment, as well as on a generalized notion of  Aubry set.

\end{abstract}
\maketitle

\begin{section}{Introduction}

Throughout the paper we deal with a Hamiltonian $H(x,p,\omega)$
  defined in $\R^N \times \R^N \times \Omega$, where
$\Omega$ is a probability space which is separable, in a suitable measure--theoretic sense.
It is assumed that $\R^N$ acts ergodically on
$\Omega$ and that $H$ satisfies a stationarity property with
respect to such action. As well known,
this frame, usually called stationary ergodic,  generalizes periodic, quasi--periodic and
almost--periodic settings. Besides the  basic convexity and growth conditions in $p$,  we require $H$ to be continuous in $x$ in
Section \ref{continuous}, and of Tonelli type in Section
\ref{Tonelli}.

The main object of our investigation   is the  stochastic Hamilton--Jacobi equation
corresponding to the critical value of $H$, say $c$,
which is in principle a random variable but stays almost surely constant, due to the ergodicity hypothesis.
It is characterized by the fact that the critical equation $H=c$
has admissible  subsolutions, but $H=a$ does not at any subcritical level, i.e.  whenever $a <c$.
By admissible subsolution we mean a random Lipschitz function that has stationary increments,
sublinear growth at infinity, and the property of being an almost everywhere subsolution in $\R^N$, almost surely with respect to $\oo$.
By convexity of the Hamiltonian, this last condition is equivalent to the notion of viscosity subsolution.


It is well known that the critical value is the unique level for
which the corresponding Hamilton--Jacobi equation may have
viscosity solutions, also termed correctors for the role they play
in associated homogenization problems. The issue of finding
characterizing conditions for the existence of correctors has been
addressed in \cite{DS12, DS11} by means of  a stochastic version
of  weak KAM Theory and of an adaptation of the metric techniques developed
for deterministic  Hamilton--Jacobi equations.  The cornerstone of
this approach is a random notion of Aubry set, defined following
the basic idea that the real invariant objects to look  at in the
analysis of stationary ergodic equations are not any more single
points of the state space but instead random closed stationary
sets.  It has been  characterized in the aforementioned references
as the minimal closed stationary random set for which there exists
an admissible subsolution strict in its complement, at least in a
weak sense, see Section \ref{section Aubry set} for precise
definitions and results.

The present paper fits into the same line of research. It
specifically aims  at  generalizing the  previous characterization
result in two directions.  On one side one would like to get
 admissible subsolutions  strict outside the random Aubry set, but in a stronger more classical sense;
on the other, assuming the Hamiltonian to be Tonelli, to show that
such random functions can be taken almost surely of class
$\CC^{1,1}$ in the whole of $\R^N$. This is actually achieved by
making use of the associated positive and negative Lax--Oleinik
semigroups, which seems  new in this context, and  shows the
effectiveness of these tools in the stationary ergodic setting,
opening  the door to a fruitful use of it for the study of other
topics in the field, as well.

 The sought  generalizations  can be, in particular,  regarded as a
  step  forward to
deduce comparison principles for the critical equation. In the
periodic case,  for instance, similar results permit to show that
the Aubry set is a uniqueness set for the critical equation,
meaning that two admissible solutions  agreeing  on it are in fact
the same. Here, as usual,  the problem gets  more involved since
it is not still clear how to exploit the weak form of compactness
encoded in the stationary ergodic model.

The  connection between regularity of subsolutions and the
property of being strict in some distinguished regions can be understood if we think of
finding regular subsolutions of the Hamilton--Jacobi equation at a supercritical level. In this
case, the problem becomes relatively easy and does not
require any special theory to be developed since, roughly
speaking,  there is some space  to make a direct regularization
starting from any globally strict subsolution, for example through
mollification, without violating the subsolution property. When
instead no strict subsolutions are available to start with,
the problem becomes difficult and requires a deeper understanding
of the setup.

To give account of the main results on the subject in the
deterministic case, we recall that  an initial crucial step has
been to realize that the obstruction in getting strict
subsolution for critical equations is not spread out indistinctly
on the whole ambient space, but is concentrated around a
specific set named after Aubry. Using this information,    the
existence of $C^1$ critical subsolutions have been first proved in
\cite{FaSic03} for Hamiltonians Lipschitz--continuous in the
state variable through a technique combining  partitions of unity
and coverings, which non surprisingly requires  quite laborious
estimates in proximity of the Aubry set.

Next, a relevant progress has been made in \cite{Be} where
$C^{1,1}$ subsolutions, which is   the optimal attainable
regularity, have been found, at least when the Hamiltonian is
Tonelli and the ambient space compact, through a simpler a more
powerful procedure based on  a double alternate application of
positive and negative  Lax--Oleinik semigroup. Within this
approach, the difficulty of dealing with the Aubry set is bypassed
thanks to the fact  that the action of such semigroups do not
affect critical subsolutions on the Aubry set. A nontrivial
extension of this result to noncompact setting has been more
recently provided in \cite{FaFiRi} by means of countable many
alternative applications of  Lax--Oleinik semigroups.

To further illustrate our results avoiding technical
complications, we assume in the remainder of the Introduction   the
critical value to be
 0 , which is not restrictive up to adding a constant to the Hamiltonian.
Our  main achievements are the following:
first, we provide a construction of the random Aubry set that
simplifies the one given in \cite{DS12} and that allows us to get
rid of a restrictive  condition therein assumed (see \eqref{C} in
Section \ref{section Aubry set}), still keeping the crucial
property  of existence of an admissible critical subsolution,
weakly strict outside it.  The crucial improvement  with respect
to the analysis performed in \cite{DS12}  being  that the main
tools used here are Lax--Oleinik semigroups instead of the
critical semidistance.

Secondly,  as already pointed out, we  in addition establish  the
existence of admissible critical subsolutions that are strict, in
the usual and stronger sense, outside the random Aubry set. More
precisely, we show that any given weakly strict critical
subsolution can be approximated, with respect to the $L^{\infty}$
norm in $\R^N$, by a critical subsolution that is strict outside
the random Aubry set, see Theorem \ref{teo strict subsol}. The key
point here is the discovery of the fact that the negative
Lax--Oleinik semigroup, when acting on a weakly strict
subsolution, produces a 1--parameter family of admissible critical
subsolutions that is {\em strictly} increasing outside the random
Aubry set, see Proposition \ref{lemma increasing}. In the end, the
sought strict subsolution is defined, as usual in the topic,
through infinite convex combination of the critical subsolutions
obtained by applying the negative Lax--Oleinik semigroup to the
initial weakly strict one at suitably small times, see Theorem
\ref{teo strict subsol}. Via a standard argument, is then easy to
prove that strict critical subsolutions are dense, with respect to
a suitable metric, among all critical admissible ones.

When the Hamiltonian is additionally assumed of Tonelli type, such strict critical subsolutions can be
 taken of class $\CC^{1,1}$ in $\R^N$. This is obtained by first deducing a further invariance of Lax--Oleinik
  semigroups with respect to strict critical subsolutions, and  then by    applying to the
 stationary ergodic environment the regularizing procedure  due to Bernard, see  \cite{Be}. The lack of
compactness of the ground space does not affect the method since,
under our assumptions,
 we have  a global control on $\R^N$ of the semiconcavity or
 semiconvexity constants of the subsolutions generated in
 the construction.  See Remark \ref{tonton2} for more comments on this issue.

The paper is  divided in three sections. Section \ref{sez basic}
has introductory character: we give notation, terminology and we
recall some basic mathematical facts that will be used throughout
the paper. Further we provide a brief presentation of  the
stationary ergodic setting and of the corresponding stochastic
Hamilton--Jacobi equations, and  present the salient properties of
the positive and negative Lax--Oleinik semigroups in random
environments.
Section 3 is about continuous Hamiltonians. In the first
subsection we give the definition of random Aubry set and we prove
the existence of an admissible critical subsolution that is weakly
strict outside it. In the second subsection we strengthen these
results by showing existence and density of strict critical
subsolutions. The final section is devoted to Hamiltonians of
Tonelli type. In the first subsection we list some additional
properties enjoyed by the Lax--Oleinik semigroups and by the random
Aubry set in our setting, while in the second one we present and
apply Bernard's method to the case at issue and we derive the
announced results about regular strict subsolutions.
\end{section}

\begin{section}{Basic material}\label{sez basic}

\begin{subsection}{Notations and preliminaries.}
We write below a list of symbols used throughout this paper.
\[
\begin{array}{ll}
N & \hbox{an integer number}\\
B_R(x_0) &\hbox{the closed ball in $\R^N$
centered at $x_0$ of radius $R$}\\
B_R & \hbox{the closed ball in $\R^N$
centered at $0$ of radius $R$}\\
\langle\,\cdot\;, \cdot\,\rangle & \hbox{the scalar product in
$\R^N$} \\
|\cdot| & \hbox{the Euclidean norm in $\R^N$}\\
\R_+& \hbox{the set of nonnegative real numbers}\\
\mathcal{B}(\R^N) & \hbox{the $\sigma$--algebra of Borel subsets
of $\R^N$}\\
\cchi_E & \hbox{the characteristic function of the set $E$}\\
\end{array}
\]

\vspace{1ex} Given a subset $U$ of $\R^N$, we denote by $\overline
U$ its closure. We furthermore say that $U$ is {\em compactly
contained} in a subset $V$ of $\R^N$ if $\overline U$ is compact
and contained in $V$. If $E$ is a Lebesgue measurable subset of
$\R^N$, we denote by $|E|$ its $N$--dimensional Lebesgue measure,
and qualify $E$ as {\em negligible} whenever $|E|=0$. We say that
a property holds {\em almost everywhere} ($a.e.$ for short) on
$\R^N$ if it holds up to a negligible set. Given a function $u$
defined in $\R^N$,  we will write $u\in C^{1,1}(B_r(x_0))$ to mean
that $u$ is of class $C^{1,1}$ on $B_r(x_0)$. We will denote by
$\D{Lip}(u;B_r(x_0))$ and $\D{Lip}(Du;B_r(x_0))$ the Lipschitz
constant of $u$  and $Du$ in $B_r(x_0)$, respectively.\par

With the term {\em curve} we will
refer to an absolutely continuous function from some given interval
$[a,b]$ to $\R^N$.
%

A $\CC^1$  function $\psi$ is
said to be {\em supertangent} (resp. {\em subtangent}) to a continuous function $u$ at a
point $x_0$ if $\psi(x_0)=u(x_0)$ and $\psi \geq u$
( resp. $\psi \leq u$) locally around $x_0$. The (possibly empty)
set made up by the differentials of supertangents (resp.
subtangents) at $x_0$ is called  superdifferential (resp.
subdifferential)  and denoted by $D^+u(x_0)$ (resp. $D^-u(x_0)$).
We recall that if both super and subdifferential are nonempty then
$u$ is differentiable at $x_0$.

Given a locally Lipschitz function $u$ in an open subset $U$ of $\R^N$, we will denote by $\partial^*u(x_0)$ the set of
{\em reachable gradients} of $u$ at $x_0$, that is the set
\[
\partial^* u(x_0)=\{\lim_n D u(x_n)\,:\,\hbox{$u$ is differentiable at $x_n$, $x_n\to x_0$}\,\},
\]
while the {\em Clarke's generalized gradient} $\partial^c u(x_0)$ is the closed convex hull of $\partial^* u(x_0)$.
The set $\partial^c u(x_0)$ contains both $D^+u(x_0)$ and $D^-u(x_0)$, in particular $Du(x_0)\in \partial^c u(x_0)$ at any differentiability point $x_0$ of $u$.
We recall that the set--valued map $x\mapsto\partial^c u(x)$ is upper semicontinuous with respect to set inclusion. When $\partial^c u(x_0)$ reduces to a singleton, the function $u$ is said to be {\em strictly differentiable} at that point. In this instance, $u$ is differentiable at $x_0$ and its gradient is continuous at $x_0$. When  $u$ depends on
a time and space variable, indicated by $t$ and $x$, respectively,
we will denote by $\partial^c_t u(x,t)$ the Clarke's generalized gradient of the function $u(\cdot,x)$ at $t$, and by
$\partial^c_x u(x,t)$ the Clarke's
generalized gradient of $u(t,\cdot)$ at $x$.
We refer the reader to \cite{Cl} for a detailed treatment of the subject.

A function $u$ will be said to be {\em semiconcave} on an open subset $U$ of $\R^N$ if for every $x\in U$ there exists  a vector $p_x\in\R^N$ such that
\begin{equation*}\label{cond semiconcavita}
 u(y)-u(x)\leqslant \langle p_x,y-x\rangle + |y-x|\,\Theta(|y-x|)\qquad\hbox{for every $y\in U$,}
\end{equation*}
where $\Theta$ is a modulus. The vectors $p_x$ satisfying such inequality are precisely the elements of $D^+ u(x)$, which is thus always nonempty in $U$. Moreover, $\partial^c u(x)=D^+u(x)$ for every $x\in U$, yielding in particular that $Du$ is continuous on its domain of definition in $U$, see \cite{CaSi}.
Finally, we say that a  function $u$ is  {\em semiconvex} if $-u$ is semiconcave.


\smallskip

Throughout the paper, $(\Omega,\F, \PP)$ will denote a {\em
separable probability space},  where $\PP$ is the probability
measure and $\F$ the $\sigma$--algebra of $\PP$--measurable sets.
Here separable is understood in the measure theoretic sense,
meaning that the Hilbert space $L^2(\Omega)$ is separable, cf.
\cite{Tsi01} also for other equivalent definitions. A property
will be said to hold {\em almost surely} ($a.s.$ for short) on
$\Omega$ if it holds up to a subset of probability 0. We will
indicate by $L^p(\Omega)$, $p\geq 1$, the usual Lebesgue space on
$\Omega$ with respect to $\PP$. If $f\in L^1(\Omega)$, we write
$\EE(f)$ for the mean of $f$ on $\Omega$, i.e. the quantity
$\int_\Omega f(\omega)\,\dd \PP(\omega)$.

We qualify  as {\em measurable} a map from $\Omega$ to itself,  or
 to a topological space $\M$ with Borel
$\sigma$--algebra $\mathcal{B}(\M)$,  if the inverse image of any
set in $\F$ or in $\mathcal{B}(\M)$  belongs to $\F$. The latter
will be also called {\em random variable} with values in $\M$.
\par

We will be interested in the case when the range of a
random variable is the {\em Polish  space} (i.e. a complete and
separable metric space) $\D C(\R^N)$ of continuous real functions on $\R^N$, endowed
with a metric $d$ inducing the topology of uniform convergence on
compact subsets of $\R^N$. It can be for instance  defined by
\begin{equation}\label{def d}
d(f,g):=\sum_{n=1}^\infty
\frac{1}{2^n}\,\frac{\|f-g\|_{L^\infty(B_n)}}{\|f-g\|_{L^\infty(B_n)}+1}\qquad
f,g\in\D C(\R^k).
\end{equation}
We will furthermore denote by  $\D{Lip}_\vartheta(\R^n)$ the subspace of $\D C(\R^N)$ made up of
Lipschitz--continuous real functions with Lipschitz constant less
than or equal to $\vartheta>0$. We will briefly say
{\em Lipschitz random function} to mean a  $\vartheta$--Lipschitz
random function for some $\vartheta>0$.

Let $(f_n)_n$ be a sequence of random variables taking values in $\D C(\R^N)$. We will say
that $f_n$ converge to $f$ {\em in probability} if, for every
$\eps>0$,
\[
\PP\left(\{\omega\in\Omega\,:\,d(f_n(\omega),f(\omega))>\eps\}\right)\to
0\quad\hbox{as $n\to +\infty.$}
\]
The limit $f$ is still a random variable. Since $\D C(\R^N)$ is a
separable metric space, almost sure convergence, i.e.
$d\left(f_n(\omega),f(\omega)\right)\rightarrow 0$ a.s. in $\omega$,
implies convergence in
probability, while the converse is not true in general. However,
the following characterization holds:

\begin{teorema}\label{tmea}
Let $f_n,f$ be random variables with values in $\D C(\R^N)$. Then $f_n\to
f$ in probability if and only if every subsequence $(f_{n_k})_k$
has  a  subsequence converging to $f$ a.s..
\end{teorema}

We denote  by $L^0(\Omega;\D C(\R^N))$ the space made up by the
equivalence classes of random variables with values in $\D C(\R^N)$  for
the relation of  almost sure equality. For every $f,g\in
L^0(\Omega;\D C(\R^N))$, we set
\[
\mu(f,g):=\inf\{\eps\geq
0\,:\,\PP\big(\{\omega\in\Omega\,:\,d(f(\omega),g(\omega))>\eps\}\big)\leq\eps\}.
\]

\begin{teorema}\label{teo Ky Fan}
 $\mu$ is a metric, named after Ky Fan, which metrizes
convergence in probability, i.e. $\mu(f_n,f)\to 0$ if and only
if $f_n\to f$ in probability, and turns $L^0(\Omega;\D C(\R^N))$ into a
Polish space.
\end{teorema}


\end{subsection}

\begin{subsection}{Stationary ergodic setting.}\label{ergo}

An {\em $N$--dimensional dynamical system} $(\tau_x)_{x\in\R^N}$
is defined as a family of mappings $\tau_x:\Omega\to\Omega$ which
satisfy the following properties:
\begin{enumerate}
\item[{\em (1)}] the {\em group property:} $\tau_0=id$,\quad
$\tau_{x+y}=\tau_x\comp\tau_y$;

\item[{\em (2)}] the mappings $\tau_x:\Omega\to\Omega$ are
measurable and measure preserving, i.e. $\PP(\tau_x E)=\PP(E)$ for
every $E\in\F$;

\item[{\em (3)}] the map $(x,\omega)\mapsto \tau_x\omega$ from
$\R^N\times\Omega$ to $\Omega$ is jointly measurable, i.e.
measurable  with respect to the product $\sigma$--algebra
$\mathcal B (\R^N)\otimes\F$.
\end{enumerate}

We will moreover assume that $(\tau_x)_{x\in\R^N}$ is {\em
ergodic,} i.e. that one of the following equivalent conditions
hold:
\begin{itemize}
\item[{\em (i)}] every measurable function $f$ defined on $\Omega$
such that, for every $x\in\R^N$, $f(\tau_x\omega)=f(\omega)$ a.s.
in $\Omega$, is almost surely constant; \item[{\em (ii)}] every
set $A\in\F$ such that $\PP(\tau_x A\,\Delta\, A)=0$ for every
$x\in\R^N$ has probability either 0 or 1, where $\Delta$ stands for
the symmetric difference.
\end{itemize}

Given   a random variable $f:\Omega\to\R$, for any fixed
$\omega\in\Omega$ the function $x\mapsto f(\tau_x\omega)$ is said
to be a {\em realization of $f$.} The following properties follow
from Fubini's Theorem, see \cite{JiKoOl}:  if $f\in L^p(\Omega)$,
then $\PP$--almost all its realizations belong to
$L^p_{loc}(\R^N)$; if $f_n \rightarrow f$ in $L^p(\Omega)$, then
$\PP$--almost all realizations of $f_n$ converge to the
corresponding realization of $f$ in $L^p_{loc}(\R^N)$. The
Lebesgue spaces on $\R^N$ are understood with respect to the
Lebesgue measure.\par

%

%
%

A jointly measurable function $v$ defined in $\R^N \times \Omega$
is said {\em stationary } if, for every $z \in \R^N$, there exists
a set $\Omega_z$ with probability $1$ such that for every
$\omega\in\Omega_z$
\[ v(\cdot + z, \omega)= v(\cdot, \tau_z \omega) \quad\text{on $\R^N$}
\]

With the term {\em (graph--measurable ) random set} we indicate  a
set--valued function $X:\Omega\to\mathcal B(\R^N)$ with
\[
\Gamma(X):=\left\{(x,\omega)\in\R^N\times\Omega\,:\,x\in
X(\omega)\,\right\}
\]
jointly measurable in $\R^N\times\Omega$. A random set $X$ will be
qualified as {\em stationary} if for every  for every $z\in\R^N$,
there exists a set $\Omega_z$ of probability 1 such that
\begin{equation}\label{def stationary set}
X(\tau_z\omega)=X(\omega)-z\qquad\hbox{for every
$\omega\in\Omega_z$.}
\end{equation}

We use a stronger notion of measurability, which is usually  named
in the literature after Effros,  to define a {\em closed random
set}, say $X(\omega)$. Namely we require $X(\omega)$ to be a
closed subset of $\R^N$ for any $\omega$ and
\begin{equation*}\label{effros}
   \{\omega \, : \, X(\omega) \cap K \neq \emptyset \} \in \F
\end{equation*}
with $K$ varying  among the compact (equivalently, open) subsets
of $\R^N$. If $X(\omega)$ is measurable in this sense then
it is also graph--measurable, see \cite{Molchanov} for more
details.

A closed random set $X$ is called stationary if it, in addition,
satisfies \eqref{def stationary set}. Note that in this event the
set $\{\omega\,:\,X(\omega)\not=\emptyset\,\}$, which is
measurable by the Effros measurability of $X$, is invariant with
respect to the group of translations $(\tau_x)_{x\in\R^N}$ by
stationarity, so it has probability either 0 or 1 by the
ergodicity assumption.

The following holds, see \cite{DS09}:

\begin{prop}\label{rando} Let $f$ be a continuous  random
function and $C$ a closed subset of $\R$. Then
\begin{equation*}\label{rando_1}
   X(\omega):= \{x \, :\, f(x,\omega)\in C\}
\end{equation*}
is a closed  random set in $\R^N$. If in addition $f$ is
stationary, then $X$ is stationary.\medskip
\end{prop}

\begin{definition}\label{def v with st increments}
A random Lipschitz function $v$ is said to have
 {\em stationary increments} if, for every
    $z\in\R^N$, there exists a set $\Omega_z$ of probability 1
    such that
\begin{equation}\label{def stationary increments}
    v(x+z,\omega)-v(y+z,\omega)=v(x,\tau_z
    \omega)-v(y,\tau_z\omega)\quad\hbox{for all  $x,y\in\R^N$}
\end{equation}
    for every $\omega\in\Omega_z$. This is equivalent to requiring
    that there exists a random variable $k(\oo)$, depending on
    $z$, for which
    \begin{equation}\label{def stationary incrementsbis}
    v(\cdot+z,\omega) =v(\cdot,\tau_z
    \omega)+ k(\oo)\quad\hbox{on $\R^N$} \qquad\hbox{for every $\oo \in \Oo_z$.}
\end{equation}

\end{definition}

Let $v$ be a Lipschitz random function with stationary increments.
For every fixed $x\in\R^N$, the random variable $Dv(x,\cdot)$ is
well defined on a set of probability 1, see \cite{DS11} for the
details. Accordingly, we can define the mean $\EE(Dv(x,\cdot))$,
which is furthermore independent of $x$ by the stationary
character of $Dv$. We are interested in the case when this mean is
zero.

\begin{definition}\label{def ammissible}
A Lipschitz random function will be called {\em admissible} if it
has stationary increments and gradient with mean 0.
\end{definition}

We state a characterizations of admissible random functions and a
result that guarantees that stationary Lipschitz random functions
are admissible.

\begin{teorema}\label{teorema mean 0}
A Lipschitz random function $v$  with stationary increments has
gradient with vanishing mean if and only if it is almost surely
sublinear at infinity, namely
\begin{equation}\label{mean1}
    \lim_{|x|\to
    +\infty}\frac{v(x,\omega)}{|x|}=0\qquad\hbox{a.s. in $\omega$.}
\end{equation}
\end{teorema}
\smallskip
%

\begin{teorema}\label{teo admissible}
Any stationary Lipschitz random function $v$ is
admissible.\medskip
\end{teorema}


\end{subsection}

\begin{subsection}{Stochastic Hamilton--Jacobi equations}\label{sez HJ}

We consider an Hamiltonian
\[
H:\R^N\times\R^N\times\Omega\to\R
\]
satisfying the following conditions:

\begin{itemize}
    \item[(H1)] the map $\omega\mapsto H(\cdot,\cdot,\omega)$
    from $\Omega$ to the Polish space $C(\R^N\times\R^N)$ is
    measurable;\smallskip
    \item[(H2)] \quad$H(x,\cdot,\omega)\ \hbox{is  convex on $\R^N$}$\ \ \qquad for every $(x,\omega)\in\R^N\times\Omega$; \smallskip
     \item[(H3)]  there
     exist two superlinear functions $\alpha,\beta:\R_+\to\R$ such
     that
     \[
     \alpha\left(|p|\right)\leq H(x,p,\omega)\leq \beta\left(|p|\right)\qquad\hbox{for all
     $(x,p,\omega)\in\R^N\times\R^N\times\Omega$;}
     \]
    \item[(H4)] \quad $H(\cdot+z,\cdot,\omega)=H(\cdot,\cdot,\tau_z\omega)$ \ \qquad for
    every $(z,\omega)\in\R^N\times\Omega$.
\end{itemize}

\begin{oss}
Condition (H3) is equivalent to saying that $H$ is superlinear and
locally bounded in $p$, uniformly with respect to $(x,\omega)$. We
deduce from (H2)
\begin{equation}\label{lippo}
  | H(x,p,\omega)-H(x,q,\omega)| \leq \widetilde L_R |p-q| \quad\text{for all
  $x$, $\omega$, and $p$, $q$ in $B_R$},
\end{equation}
where $\widetilde L_R=\sup\{\, |H(x,p,\omega)|\,:\,(x,\omega) \in
\R^N \times \Omega, \,|p|\leq R+2\,\},$ which is finite thanks to
(H3).
\end{oss}
%

For every $a\in\R$, we are interested  in the stochastic
Hamilton--Jacobi equation
\begin{equation}\label{eq HJa}
H(x,Dv(x,\omega),\omega)=a\qquad\hbox{in $\R^N$.}
\end{equation}
The material we are about to expose has been already presented in
\cite{DS12, DS11, DS09}, to which we refer for the details. Here we just
recall the main items.
\par

We say that a Lipschitz random function is a {\em solution} (resp.
{\em subsolution}) of \abra{eq HJa} if it is a viscosity solution
(resp. a.e. subsolution) of \abra{eq HJa} for every $\omega$ in a
set of probability 1. We recall that, due to the convexity
assumption (H2), the notion of almost everywhere subsolution is
equivalent to the one in the viscosity sense. We refer the reader to
\cite{BCD97, Ba94} for more details on the theory of viscosity
solutions in the deterministic case. Notice that any such
subsolution is almost surely in $\D{Lip}_{\kappa_a}(\R^n)$, where
\begin{equation}\label{def kappa_a}
\kappa_a:=\sup\{\,|p|\,:\,H(x,p,\omega)\leq a\ \hbox{for some
$(x,\omega)\in\R^N\times\Omega$}\,\},
\end{equation}
which is finite thanks to (H3). We are interested in the class of
{\em admissible subsolutions}, hereafter denoted by $\S_a$, i.e.
random functions taking values in $\D{Lip}_{\kappa_a}(\R)$ with
stationary increments and zero mean gradient that are subsolutions
of \eqref{eq HJa}. An admissible solution will be also named {\em
exact corrector}, remembering its role in homogenization.

%

\medskip

We proceed by defining the {\em free} and the {\em stationary random
critical value}, denoted by $c_f(\omega)$ and $c$ respectively, as
follows:
\begin{eqnarray}
 c_{f}(\omega)&=& \inf\left\{ a\in\R\,:\, \text{\abra{eq HJa} has a subsolution
$v\in\D{Lip}(\R^N)$} \right\}, \label{def cf}\\
c &=& \inf\{a\in\R\,:\,\S_a\not =\emptyset\,\}.  \label{def c}
\end{eqnarray}
We emphasize that in definition \abra{def cf} we are considering
{\em deterministic} a.e. subsolutions $v$ of the equation \abra{eq
HJa}, where $\omega$ is treated as a fixed parameter. Furthermore,
we note that $c_f(\tau_z\omega)=c_f (\omega)$ for every
$(z,\omega)\in\R^N\times\Omega$, so that, by ergodicity, the
random variable $c_f(\omega)$ is almost surely equal to a
constant, still denoted by  $c_f$. Hereafter we will write
$\Omega_f$ for the set of probability $1$ where $c_f(\omega)$
equals $c_f$.

It is apparent by the definitions that $c \geq c_f$. The relation
of these two values with the effective Hamiltonian obtained via
the homogenization \cite{ReTa00, Souga99} is discussed in
\cite{LiSou03, DS11}.


In the sequel, we will focus our attention on the {\em critical
equation}
\begin{equation}\label{eq critica}
H(x,Dv(x,\omega),\omega)=c\qquad\hbox{in $\R^N$.}
\end{equation}
The following result holds, see \cite{LiSou03, DS11}:

\begin{teorema}
There exist admissible critical subsolutions, i.e.
$\S_c\not=\emptyset$.
\end{teorema}

Moreover, the critical equation is the only equation of the kind
\eqref{eq HJa} for which exact correctors may exist, see
\cite{DS09}.
\smallskip

We recall the main items of the so called metric method for
Hamilton--Jacobi equations, as adapted in \cite{DS11} to the
case at issue.  In next formulae we assume that $a \geq c_f$ and
$\omega \in \Omega_f$. Let
\[
Z_a(x,\omega):=\{p\,:\,H(x,p,\omega)\leq a\,\},
\]
be the $a$--sublevels of the Hamiltonian and
    \begin{equation*}\label{sigma}
    \sigma_a(x,q,\omega):=\sup\left\{\langle q,p\rangle\,:\,p \in Z_a(x,\omega)\,\right\}
    \end{equation*}
the related support functions. It comes from \abra{lippo} (cf.
Lemma 4.6 in \cite{DS09}) that, given $b>a$, we can find
$\rho=\rho(b,a)>0$ with
\begin{equation}\label{spessore}
Z_a(x,\omega)+B_\rho\subseteq Z_b(x,\omega) \qquad\hbox{for every
$(x,\omega)\in\R^N\times\Omega_f$.}
\end{equation}
It is straightforward to check that $\sigma_a$ is convex in $q$,
upper semicontinuous in $x$ and, in addition, continuous whenever
$Z_a(x,\omega)$ has nonempty interior or reduces to a point. We
extend the definition of $\sigma_a$ to
$\R^N\times\R^N\times\Omega$ by setting
$\sigma_a(\cdot,\cdot,\omega)\equiv 0$ for every
$\omega\in\Omega\setminus\Omega_f$. With this choice, the function
$\sigma_a$ is jointly measurable in $\R^N \times\R^N\times \Omega$
and enjoys the stationarity property
\[
\sigma_a(\cdot+z,\cdot,\omega)=\sigma_a(\cdot,\cdot,\tau_z\omega)\quad\hbox{for
every $z\in\R^N$ and $\omega\in\Omega$.}
\]
We define the semidistance $S_a$  as
    \begin{equation}\label{eq S}
    S_a(x,y,\omega)=\inf\left\{\int_0^1 \sigma_a(\gamma(s),\dot\gamma(s),\omega)\,\dd
    s\,:\, \gamma\in\D{Lip}_{x,y}([0,1];\R^N)\, \right\}.
    \end{equation}
The following holds, cf. \cite{FaSic03}:

\begin{prop}\label{prop S} Let $ a \geq c_f$ and $\omega\in\Omega_f$. We have:
 \begin{itemize}
     \item[\em(i)] For any $y\in\R^N$, the function $S_a(y,\cdot,\omega)$ is a subsolution of \abra{eq HJa} in $\R^N$, and also a solution in $\R^N\setminus\{y\}$.
     \item[\em(ii)] A continuous function $\phi$ is a subsolution of \abra{eq HJa} if and
     only if
 \[
 \phi(x)-\phi(y)\leq S_a(y,x,\omega)\qquad\hbox{for all
 $x,y\in\R^N$.}
 \]
 \end{itemize}
 \end{prop}
\medskip
\end{subsection}

\begin{subsection}{Positive and negative Lax--Oleinik semigroups.} The {\em Lagrangian} associated with $H$ by duality is the
function defined as
\[
L(x,q,\omega)
    :=
   \max_{p\in\R^N}\left\{\langle p,q\rangle -
   H(x,p,\omega)\right\}\quad\hbox{for all $(x,q,\omega)\in\R^N\times\R^N\times\Omega$.}
\]
As well known, the Lagrangian satisfies properties analogous to
(H1)--(H4).\smallskip

For every $t>0$ and $\omega\in\Omega$, we define a function $h_t$
on $\R^N\times\R^N$ as
\[
h_t(x,y,\omega):=\inf \left\{\int_{0}^t
L(\gamma,\dot\gamma,\omega)\,\dd
s\,:\,\gamma(0)=x,\,\gamma(t)=y\,\right\}.
\]
Curves that realize the above infimum are called {\em Lagrangian
minimizers.} It is well known that such minimizers always exist
and they are Lipschitz continuous, see \cite{BuGiHi98}.

It is known, see for instance \cite{Fathi}, that
\begin{equation}\label{eq ht>S}
S_a(y,x,\omega)=\inf_{t>0} \Big(h_t(y,x,\omega)+at\Big)
\end{equation}
for every $a\in\R$ and $\omega\in\Omega_f$. In particular, when
$\omega\in\Omega_f$ and $a\geq c_f$, we derive, in view of
Proposition \ref{prop S}, that a function $v$ is a subsolution of
\eqref{eq HJa} if and only if
\[
v(x)-v(y)\leq h_t(x,y)+a\,t\qquad\hbox{for every $x,y\in\R^N$ and
$t>0$.}
\]

For every $t>0$ and every $\kappa$--Lipschitz random function $u$,
we next define the {\em negative} and {\em positive Lax--Oleinik
semigroups} $T_t^-u$, $T_t^+u$ as follows:
\begin{eqnarray*}
(T_t^-u)(x,\omega)&:=&\inf\left\{u(y,\omega)+h_t(y,x,\omega)\,:\,y\in\R^N\,\right\}\\
(T_t^+u)(x,\omega)&:=&\sup\left\{u(y,\omega)-h_t(x,y,\omega)\,:\,y\in\R^N\,\right\}
\end{eqnarray*}
for every $(x,\omega)\in\R^N\times\Omega$.\smallskip

\begin{oss}\label{oss H check}
It is worth noticing that $T_t^+u=-\check T_t^-(-u)$, where we
have denoted by $\check T_t^-$ the negative Lax--Oleinik semigroup
associated with the Hamiltonian $\check H$ defined by
\[
\check H(x,p,\omega):=H(x,-p,\omega)\qquad\hbox{for every
$(x,p,\omega)\in\R^N\times\R^N\times\Omega$.}
\]
\end{oss}

\smallskip
We list a series of properties enjoyed by such semigroups.

\begin{prop}\label{prop B4}
Let $\vartheta>0$. Then there exists $R(\vartheta)>0$ such that
for every $\vartheta$--Lipschitz random function $u$ and every
$\omega\in\Omega$ the following properties hold:
\begin{itemize}
    \item[{\em (i)}]\quad $(T_t^-u)(x,\omega)
<\inf\left\{u(y,\omega)+h_t(y,x,\omega)\,:\,|x-y|\geq
t\,R(\vartheta)\,\right\},$\smallskip
    \item[] \quad
$(T_t^+u)(x,\omega)>\sup\left\{u(y,\omega)-h_t(y,x,\omega)\,:\,|x-y|\geq
t\,R(\vartheta)\,\right\}$;\medskip
    \item[{\em (ii)}]  the maps
$(t,x)\mapsto T_t^{\pm} u(x)$ are $R(\vartheta)$--Lipschitz
continuous on $[0,+\infty)\times\R^N$. In particular
\[
\|T_t^\pm u(\cdot,\omega)-u(\cdot,\omega)\|_{\infty}\leq
t\,R(\vartheta)\qquad\hbox{for every $t>0$;}
\]
\item[{\em (iii)}] if $u$ has stationary increments (resp. is
stationary), then $T^\pm_t u$ have stationary increments (resp.
are stationary) for any $t>0$.
\end{itemize}
\end{prop}

\smallskip

Items {\em (i), (ii)} are the ergodic stationary version of well known
estimates holding in deterministic case. Bounding function
$R(\cdot)$ solely depends on the $\alpha(\cdot)$, $\beta(\cdot)$
appearing in (H3), which are invariant with respect  to $\oo$.
Item {\em (iii)} is straightforward.

\vspace{1ex}

\begin{prop}\label{prop Tm solution}
Let $u_0$ be a Lipschitz random function. Then, for every fixed
$\omega\in\Omega$, the function $u(t,x)=(\Tm u_0)(x,\omega)$ is
the unique Lipschitz continuous viscosity solution of the time--dependent equation
\begin{eqnarray}\label{Cauchy problem}
\partial_t u+H(x,Du,\omega)=0& \hbox{in $(0,+\infty)\times \R^N$}
\end{eqnarray}
satisfying \  $u(0,x)=u_0(x,\omega)\quad\hbox{on $\R^N$.}$
\end{prop}

\vspace{2ex}

\begin{prop}\label{prop critical sol} Let $u$ be a continuous function on $\R^N$ and let $\omega$ be fixed.
The following facts hold:
\begin{itemize}
    \item[{\em (i)}] $u$ is a subsolution of \eqref{eq HJa} if and only
if $t\mapsto \Tm u(\cdot,\omega)+a\,t$ is non--decreasing;\smallskip
    \item[{\em (ii)}] $u$ is a solution of \eqref{eq HJa} if and only
if $u\equiv \Tm u(\cdot,\omega)+a\,t$ for every $t>0$.
\end{itemize}
\end{prop}

As a consequence of the previous results, we deduce the following
fact:

\begin{prop}\label{prop B subsol}
Let $u\in\S_a$. Then both $\Tp u$ and $\Tm u$ belong to $\S_a$,
for every $t>0$.
\end{prop}
\medskip
%
%
%
%
\end{subsection}
\end{section}

\begin{section}{Continuous Hamiltonians}\label{continuous}

In this section we assume $H$ to satisfy assumptions (H1)--(H4).
In the first subsection we give the definition of random Aubry set
and we prove the existence of a critical admissible subsolution
that is {\em weakly strict}, in a sense that will be clarified,
outside it. This can be regarded as a completion of \cite{DS12},
where analogous results were obtained under an additional
assumption, see \eqref{C} below. In the second subsection we show
how a strict admissible critical subsolution can be produced
starting from a weakly strict one. For the outcomes of both
subsections, we essentially use  the (negative) Lax--Oleinik
semigroup and its  properties.

From now on we will assume the stationary critical value
$c=0$, which is not restrictive up to replacing $H$ with $H-c$,
and focus our analysis on such critical level. To ease notation,
we will omit the subscript $c$ when no ambiguity is possible. In
particular, we will simply write $\S$, $S$, $\kappa$
and $\sigma$ in place of $\S_c$,  $S_c$, $\kappa_c$
and $\sigma_c$.

\begin{subsection}{Aubry set and weakly strict subsolutions}\label{section Aubry set}
We start by recalling the definition of the classical Aubry set.
\begin{definition}
The {\em classical Aubry set} $\A_f(\omega)$ is the closed stationary random set  defined as follows:
\begin{equation*}\label{aub}
    \A_f(\omega)=\{\,y\in\R^N\,:\ \liminf_{t\to +\infty}
    h_t(y,y,\omega)=0\,\}\quad\hbox{for any $\omega\in\Omega$.}
\end{equation*}
\end{definition}

The random set  $\A_f(\omega)$ is always almost surely empty when the random critical value is strictly greater than $c_f$, i.e. when $c_f<0$. When $c_f=0$, $\A_f(\omega)$ may be either almost surely nonempty or almost surely empty. This latter instance is not specific of the stationary ergodic setting: it can occur even in the periodic case.    

We proceed by presenting the notion of random Aubry set as introduced in \cite{DS12}. For this, we need to recall some facts about 
Lax formulae in the stationary ergodic
setting, see \cite{DS11, DS09} for more details. For any given almost surely nonempty stationary closed
random set $C(\omega)$ in $\R^N$ and any critical subsolution
$g\in\S$, the Lax formula is given by
\begin{equation}\label{def u}
u(x,\omega):= \inf\{g(y,\omega)+S(y,x,\omega)\,:\,y\in
C(\omega)\,\}\quad \hbox{$x\in\R^N$,}
\end{equation}
and provides another  admissible critical subsolution enjoying
some additional properties, see Proposition 4.1 in \cite{DS12}.
More precisely, $u(\cdot,\omega)$ is the maximal critical
subsolution agreeing with $g(\cdot,\omega)$ on $C(\omega)$ a.s. in
$\omega$, and, as a consequence
\begin{equation}\label{eq maximal subsol}
 \hbox{$u(\cdot,\omega)$\quad is a critical solution in $\R^N\setminus C(\omega)$\qquad a.s. in $\omega$ .}
\end{equation}
In the above formula we agree that $u(\cdot,\omega)\equiv 0$ when either
$C(\omega)=\emptyset$ or the infimum is equal to $-\infty$.

Inspired by the  periodic model case, we define the stationary
random Aubry set as follows:
\begin{definition}\label{def Aubry} A stationary and almost surely nonempty closed random set $C(\omega)$  will be called
{\em random Aubry set} and will be denoted by $\A(\omega)$ if 
\begin{itemize}
    \item[\em (i)] the extension of any  admissible critical subsolution
from  $C(\omega)$ via the Lax formula \eqref{def u} yields an
admissible critical solution;
    \item[\em (ii)] any  stationary and almost surely nonempty closed random set that satisfies the
previous property  is almost surely contained in $C(\omega)$.\smallskip
\end{itemize}
If there are no stationary and almost surely nonempty closed random sets satisfying {(i)}, we agree that the {\em Aubry set} is almost surely empty.
\end{definition}

\smallskip
We now proceed to show that the above definition is well posed.
For this, we will take on the burden of checking that, if there is some nonempty 
stationary closed random set satisfying item {\em(i)} above, then there is also a maximal one. 
This task is by no means trivial, the main difficulty being the following: the random sets needed for the construction cannot be defined by treating $\omega$ as a fixed parameter because this procedure would lead to objects that do not satisfy the proper measurability hypotheses. 

We also underline that the case of an almost surely empty random Aubry set (i.e. the case of nonexistence of nonempty 
stationary closed random set satisfying item {\em(i)} above) can actually occur, 
see for instance Example 4.10 in \cite{DS11}. 

The fact that Definition \ref{def Aubry} is well--posed has been
already shown in \cite{DS12} under the following additional
assumption:
\begin{equation}\label{C}\tag{A}
\hbox{either\quad $c_f<0$\quad or\quad $c_f=0$\  \ and\ \
$\A_f(\omega)= \emptyset$ a.s. in $\omega$.}
\end{equation}
As already announced, the novelty here is that we  get rid of \eqref{C} by
using the negative Lax--Oleinik semigroup.  

In order to state the result we aim at, we need the following 

\begin{definition}\label{def weakly strict}
Let $C(\omega)$ a stationary closed random set. A critical subsolution $v\in\S$ is said to be {\em weakly strict} in
$\R^N\setminus C(\omega)$ if the following property holds a.s. in $\omega$:
\begin{equation}\label{ineq weakly strict}
  v(x,\omega)-v(y,\omega)<S(y,x,\omega)\qquad\hbox{for every $x, y\in \R^N\setminus C(\omega)$ with $x\not =y$.}
\end{equation}
\end{definition}

Notice that the inequality \eqref{ineq weakly strict} keeps holding even if one (but only one)
between the points $x,y$ belongs to $C(\omega)$. This can be deduced from the fact
 that $S(\cdot,\cdot,\omega)$ is
a geodesic--type semidistance.
\smallskip

\begin{teorema}\label{teo Aubry set}
The definition of random Aubry set $\A(\omega)$ is well posed. Moreover, there is a
critical subsolution which is weakly strict in
$\R^N\setminus\A(\omega)$.
\end{teorema}

\medskip

Under assumption \eqref{C}, it has been provided in \cite{DS12} a
sort of dual characterization of Aubry  set as  the minimal
stationary closed random set for which there exists a critical
admissible subsolution  weakly strict in its complement. To prove
the above theorem, we try to generalize this approach. We thus
define, for each $v\in\S$, 
\begin{equation}\label{def X_u}
X_v(\omega)=\{x\in\R^N\,:\,(\Tm
v)(x,\omega)=v(x,\omega)\quad\hbox{for some $t>0$}\,\}.
\end{equation}
At first sight these objects are  not so  appealing, from a
mathematical point of view,  since it does not seem possible to
determine their randomness properties or even the topological
features of their images, as $\oo$ varies in $\Oo$. Nevertheless,
we keep under control our uneasiness and go on constructing,
through suitable infinite convex combinations, a distinguished
critical subsolution, denoted by $w$, for which it is possible to define a 
stationary closed random set almost surely agreeing with $X_w(\omega)$.

We proceed by showing a property enjoyed by such sets $X_v$ that will be relevant for our analysis. 
The following fact will be exploited in the proof and later on in the paper: when $v(\cdot,\omega)$
is a critical subsolution, the equality $(\Tm
v)(x,\omega)=v(x,\omega)$ for some $t>0$ implies the existence of a
point $y\in\R^N$ satisfying
\[
h_t(y,x,\omega)=v(x,\omega)-v(y,\omega)=S(y,x,\omega).
\]
The last equality follows from Proposition \ref{prop S} and
\eqref{eq ht>S} with $a=0$.

\begin{prop}\label{prop supersol on Xv}
Let $v\in\S$. Then there exists a set $\Omega_v$ of probability 1 such that, for every $\omega\in\Omega_v$ and  $x_0 \in X_v(\omega)$,  any $\CC^1$  subtangent $\psi$ to
$v(\cdot,\omega)$ at $x_0$  satisfies $H(x_0,D\psi(x_0),\oo) \geq
0$.
\end{prop}

\begin{dimo}
Let us fix $\omega$ in a set $\widehat\Omega$ of probability 1
such that $v(\cdot,\omega)$ is a critical subsolution and pick
$x_0\in X_v(\omega)$. By definition of $X_v(\omega)$, we find
$t>0$ and $y\in\R^N$ with
\begin{equation}\label{Xv1}
    v(x_0,\omega)=v(y,\omega)+h_t(y,x_0,\omega)=v(y,\omega)+S(y,x_0,\omega).
\end{equation}
This implies that the function $\overline
v(\cdot):=S(y,\cdot,\omega)+v(y,\omega)$ touches $v(\cdot,\omega)$
from above at $x_0$, so any subtangent to $v(\cdot,\omega)$ at
$x_0$ is also a subtangent to $\overline v$ at that point.

If $y\not=x_0$, the assertion follows since $\overline v$ is a
critical solution in $\R^N \setminus \{y\}$. If instead  $y=x_0$,
then the we deduce from \eqref{Xv1} the equality
$h_t(x_0,x_0,\omega)=S(x_0,x_0,\omega)=0$. This means that there exists a closed curve 
$\gamma$  parameterized in the interval $[0,t]$ and with base point at $x_0$ along which the action of $L(\cdot,\cdot,\omega)$ 
is zero. By moving along $\gamma$ $k$--times we infer that $h_{k\,t}
(x_0,x_0,\omega)=0$ for any $k \in \N$. This, in turn, implies
\[ 
\liminf_{s \to + \infty} h_s(x_0,x_0,\omega)=0,
\]
and this relation  is possible only  when $c_f=0$ and
$x_0\in\A_f(\omega)$.  In this instance $\overline v$ is a
critical solution on the whole $\R^N$, see \cite{Fathi} (or
\cite{DavZav} for  a proof in the   case of a convex Hamiltonian
just continuous). This concludes the proof.
\end{dimo}

\medskip

Given $v\in\S$, we set
\begin{equation}\label{hatv}
  \hat v(x,\omega):=v(x,\omega)-v(0,\omega)\qquad\hbox{for every $(x,\omega)\in\R^N\times\Omega$.}
\end{equation}
The  set made up by  critical subsolutions obtained in this way,
namely
\[
\widehat\S:=\{\hat v\in\S\,:\,\hat v(0,\omega)=0\ \hbox{for every
$\omega$}\}
\]
is a subspace of $L^0(\Omega,\D{C}(\R^N))$, in particular it is
separable with respect to the Ky Fan metric by Theorem \ref{teo Ky
Fan}. There thus  exists a sequence of Lipschitz random functions
$(v_n)_n$  dense in $\widehat S$ with respect to convergence in
probability, which implies, in view of Theorem \ref{tmea}, that it
is also dense for the almost sure convergence in $\D C(\R^N)$. We
set
\begin{equation}\label{def w}
w(x,\omega)=\sum_{n=1}^{+\infty}\frac{1}{2^n} v_n(x,\omega)\qquad\hbox
{for every $(x,\omega)\in\R^N\times\Omega$.}
\end{equation}

\medskip

The next result illustrates  the first crucial  property enjoyed
by $w$.
\begin{prop}\label{prop1 w}
For every $v\in\S$ there exists a set $\Omega_v$ of probability 1 such that
\[
X_w(\omega)\subseteq X_v(\omega)\qquad\hbox{for every $\omega\in\Omega_v$.}
\]
\end{prop}
\begin{dimo}
 The proof is divided in two parts. First we prove the assertion
 for any element of the sequence $(v_n)_n$,
 then, in the second step, we extend it, by density,  to  all random
  functions in $\widehat\S$. This is actually enough, for
 $X_v(\oo)= X_{\hat v}(\oo)$, for any $\oo$, whenever $v$ and
$\hat v$ are in the relation given by \eqref{hatv}. We set
\[ \widehat\Omega:=\{\omega\in\Omega\,:\,v_n(\cdot,\omega)\ \hbox{is
a critical subsolution for every $n\in\N$}\,\}\]  and pick $\oo
\in \widehat\Omega$,  $x\in X_w(\omega)$. Then there exists $t>0$
such that $(\Tm w)(x,\omega)=w(x,\omega)$. By the very definition
of $\Tm$
\[
(\Tm w)(x,\omega)\geq \sum_{n=1}^{+\infty}\frac{1}{2^n}\,(\Tm
v_n)(x,\omega),\] and, combining this information with the
monotonic character  of the action of the negative Lax-Oleinik
semigroup on critical subsolutions pointed out in Proposition
\ref{prop B subsol}, we get\[ (\Tm w)(x,\omega)\geq
\sum_{n=1}^{+\infty}\frac{1}{2^n}\,(\Tm v_n)(x,\omega)\geq
\sum_{n=1}^{+\infty}\frac{1}{2^n}\, v_n(x,\omega)=w(x,\omega),
\]
discovering in the end  that all the inequalities in the above
formula have actually to be equalities.  In particular we infer
\begin{equation}\label{eq X_vn}
(\Tm v_n)(x,\omega)=v_n(x,\omega)\qquad\hbox{for every $n\in\N$,}
\end{equation}
which  proves that $X_w(\omega)\subseteq X_{v_n}(\omega)$\quad for
every $\omega\in\widehat\Omega$ and $n\in\N$. Now let
$v\in\widehat\S$. Then there exist a subsequence $(v_{n_k})_k$ such that 
\begin{equation}\label{conve}
   v_{n_k}(\cdot,\omega)\ucv v(\cdot,\omega)
\end{equation}
for every $\omega$ in a set
$\Omega_v\subseteq\widehat\Omega$ of probability 1. Given
$\omega\in\Omega_v$ and a point $x\in X_w(\omega)$, we have, in
force of the first part of the proof, that \eqref{eq X_vn} holds
for some $t >0$, from which we  derive 
$\Tm v(x,\omega)=v(x,\omega)$ in view of \eqref{conve}. This shows that $x\in X_v(\omega)$, as
desired.
\end{dimo}

\medskip
\begin{oss}\label{remw} The
argument in the previous  proof  actually allows to establish the following more general result: 
for any given $v \in \S$ and  $t>0$
\[X_v^t(\oo) \subset X_w^t(\oo) \qquad\hbox{for a.s. $\oo$,}\]
where
\begin{equation}\label{ww}
    X_v^t(\oo) = \{x\in\R^N\,:\,(\Tm v)(x,\omega)=v(x,\omega)\}
\end{equation}
and $ X_w^t(\oo)$ is defined similarly.
\end{oss}

\medskip
In the next proposition we show that the action of the negative Lax--Oleinik semigroup does not affect
the values of $w$ on $X_w(\oo)$ not only for some positive $t$, as required
in the original definition of such a  set, see \eqref{def X_u},
but actually for {\em all} $t>0$, at least a.s. in $\oo$. 
\smallskip

\begin{prop}\label{prop2 w}
\[
X_w(\omega)=\{x\in\R^N\,:\,(\Tm
w)(x,\omega)=w(x,\omega) \;\hbox{for every $t>0$}\}\quad\hbox{a.s. in
$\omega$.}
\]
\end{prop}

\begin{dimo}
Let $\widehat\Omega$ be a set of probability 1 such that $w(\cdot,\omega)$ is a critical subsolution. According to the monotonicity  property stated in Proposition \ref{prop critical sol}, the following implication holds for every $\omega\in\widehat\Omega$:
\begin{equation}\label{2ww00}
   (\Tm w)(x,\oo)=w(x,\oo) \Rightarrow T^-_s
w(x,\oo)=w(x,\oo) \quad\hbox{ for $ s \in [0,t]$.}
\end{equation}
To prove the assertion, it will be enough to show that, for every $\omega$ in a set of probability 1, 
\[ 
\sup\{t >0\,:\, \Tm w(x,\oo)=w(x,\oo)\} = + \infty\qquad\hbox{for every $x \in
X_w(\omega)$.}
\]
To this aim, we consider a sequence $(t_n)_n$ dense
in $[0,+ \infty)$, and set  $v_n =T^-_{t_n} w$.  According to
Remark \ref{remw}, there is a set  $\Omega'\subseteq\widehat \Oo$ of probability 1
with
\begin{equation}\label{2ww1}
    X^{t_n}_{w}(\oo) \subset X^{t_n}_{v_n}(\oo) \qquad \hbox{for
any $n \in  \N$ and  $\oo \in \Omega'$},
\end{equation}
see \eqref{ww} for the notation.  Fix $\oo \in \widehat\Oo$, $x
\in X_w(\oo)$, and assume for purposes of  contradiction
\[\max\{t >0\,:\, \Tm w(x,\oo)=w(x,\oo)\} =: s  < \infty,\]
 then select   $t_k$ with
 \begin{equation}\label{2ww11}
     t_k < s < 2\,t_k.
\end{equation}
 By \eqref{2ww00}
 \begin{equation}\label{2ww2}
    (T^-_{t_k}w)(x,\oo)= w(x,\oo),
\end{equation}
or, in other terms, $x \in  X^{t_k}_{w}(\oo)$. In view of 
\eqref{2ww1} we get $x \in X^{t_k}_{v_k}(\oo)$, i.e.
\[ (T^-_{t_k+t_k}w)(x,\oo) = (T^-_{t_k}w)(x,\oo).\]
Combining it with \eqref{2ww2}, we get
\[(T^-_{2 t_k}w)(x,\oo) = w(x,\oo),\]
which is in contrast with \eqref{2ww11} and the maximality
property of $s$. 
\end{dimo}

\medskip

We go on gathering some more information on $X_w$.

\begin{prop}\label{prop Af}
Let $c_f=0$. There exists a set $\widehat \Omega$ of probability 1
such that
$$\A_f(\omega)\subseteq X_w(\omega)\qquad\hbox{for every $\omega\in\widehat\Omega$.}$$
\end{prop}

\begin{dimo}
Let $\widehat\Omega$ be a set of probability 1 such that
$w(\cdot,\omega)$ is a critical subsolution. Then,  exploiting
Proposition \ref{prop critical sol} and the very definition of
$\Tm w$ , we have, for every $t>0$ and $x \in \R^N$,
\begin{equation}\label{Af1}
    w(x,\omega)\leq (\Tm w)(x,\omega)\leq \liminf_{t\to +\infty}
w(x,\omega)+h_t(x,x,\omega).
\end{equation}
When  $c_f=0$ and $x\in\A_f(\omega)$ 
\begin{equation}\label{Af2}
  \liminf_{t\to +\infty} h_t(x,x,\omega)=0,
\end{equation}
so that, combining \eqref{Af1} and \eqref{Af2}, we get in the end
$w(x,\omega)= (\Tm w)(x,\omega)$.
\end{dimo}

\bigskip

\medskip

We pause our  analysis on $w$ to  derive a general
characterization of weakly strict subsolutions in terms of strict
monotonicity  of the action of the negative Lax--Oleinik
semigroup.

\medskip

\begin{lemma}\label{lemma pre increasing} Let $v
\in \S$ and $C(\omega)$ be a stationary closed random set. If $c_f=0$, we additionally assume
that $\A_f(\omega)\subseteq C(\omega)$ a.s. in
$\omega$. Then $v$ is weakly strict in $\R^N\setminus
C(\omega)$ if and only if there exists a set $\widehat{\Oo}$ of
probability $1$  such that
\begin{equation}\label{preincr1}
    (T^-_t w)(y,\omega)> w(y,\omega) \quad\text{for any $y \in \R^N \setminus
    C(\oo)$, $\oo \in \widehat{\Oo}$  and $t>0$.}
\end{equation}
\end{lemma}

\begin{dimo}
Let  $\widehat\Oo$ be a set of probability 1 made up of elements
$\omega\in\Omega$ for which $w(\cdot,\oo)$ is a critical
subsolution and enjoys \eqref{ineq weakly strict}. Throughout the
proof, $\omega$ will denote a fixed element of  $\widehat\Oo$.
Were \eqref{preincr1}  not true, we should have by Proposition
\ref{prop critical sol}
\[(T^-_t w)(y,\omega)=w(y,\omega),\]
for some $t >0$ and $y \not \in C(\oo)$, and by the principle we
have written down  before the statement of Proposition \ref{prop supersol on Xv}, this should in turn imply
\begin{equation}\label{incre1}
    w(y,\omega)-w(z,\omega)=S(z,y,\omega) \qquad\hbox{for some $z \in
\R^N$,}
\end{equation}
in contrast with the fact that $w(\cdot,\omega)$ is weakly strict
and $y\not\in C(\omega)$.

The converse implication will be also proved by contradiction.
Assume that \eqref{incre1} holds for some $y,\,z$ not belonging
to $C(\oo)$.  Exploiting that $S(\cdot,\cdot,\omega)$ is a
geodesic semidistance, we can assume, up to moving  $y$, that
there is a curve $\gamma:[0,1] \to \R^N$ joining $z$ to $y$ and
with support disjoint from $C(\oo)$ such that
\[
w(y,\oo) - w(z,\oo)= \int_0^1 \sigma (\gamma,\dot\gamma,\oo) \, \dd
s.
\] 
Since the support  of $\gamma$ is also by assumption disjoint
from $\A_f(\oo)$ in the case that $c_f=0$, we can provide such curve a
$0$--Lagrangian change of parameter in a bounded interval $[0,t]$,
for some $t >0$, see \cite{D1-06}, \cite{DS06}. Namely, we can
determine  a reparametrization $\xi$ with
\[\sigma(\xi(s), \dot\xi(s),\oo) =L(\xi(s), \dot\xi(s),\oo) \qquad\hbox{for
a.e. $s \in [0,t]$.}\] We deduce
\[w(y,\oo) = w(z,\oo) + \int_0^t  L(\gamma,\dot\gamma,\oo) \, \dd
s \geq (\Tm w)(y,\oo),\] but since the opposite inequality comes from
Proposition \ref{prop critical sol}, the above formula must actually
hold with equality, yielding a contradiction.

\end{dimo}

\medskip

We are now in position to give the \\

\noindent{\bf Proof of Theorem \ref{teo Aubry set}.}\ Let $w$ be defined via \eqref{def w} and set 
\[
C(\omega)=\bigcap_{n\in\N}\{x\in\R^N\,:\,(T^{-}_{s_n} w)(x,\omega)=w(x,\omega)\,\}\quad\hbox{for every $\omega$,}
\]
where $(s_n)_n$ is any dense sequence in $\R_+$. It is easily seen that $C(\omega)$ is stationary closed random set. 
Moreover, by the monotonicity property stated in Proposition \ref{prop critical sol} and by Proposition \ref{prop2 w},
\begin{equation}\label{eq prop2 w}
 C(\omega)=X_w(\omega)\qquad\hbox{a.s. in $\omega$.}
\end{equation}
We claim that $w$ is weakly strict in $\R^N\setminus C(\omega)$ and that $C(\omega)$ is the Aubry set.

Let us start with the weakly strict character of $w$ in $\R^N\setminus C(\omega)$. In view of \eqref{eq prop2 w} and of Propositions \ref{prop Af}, for every $\omega$ in a set of probability 1 we have  
\begin{equation}\label{ineq strict}
 (\Tm w)(y,\omega)>w(y,\omega)\qquad\hbox{for every $y\in \R^N\setminus C(\omega)$ and $t>0$}
\end{equation}
and $\A_f(\omega)\subseteq C(\omega)$ if $c_f=0$.
By Lemma \ref{lemma pre increasing}, we conclude that $w$ is weakly strict in $\R^N\setminus C(\omega)$.
 
Let us now show that $C(\omega)$ is the Aubry set, i.e. it is the maximal stationary closed random set satisfying item {\em (i)} in Definition \ref{def Aubry}. We start by proving the maximality property. 
Let $\widetilde C(\omega)$ be a nonempty stationary closed random set satisfying item {\em (i)} in Definition \ref{def Aubry} and set 
\[
u(x,\omega):= \inf\{w(y,\omega)+S(y,x,\omega)\,:\,y\in
\widetilde C(\omega)\,\},\qquad \hbox{$x\in\R^N$.}
\]
Since $u$ is an admissible critical solution, we infer from Proposition \ref{prop critical sol} that the following equality holds for every $\omega$ in a set of probability 1:
\begin{equation}\label{Au1}
    \Tm u(\cdot,\omega) \equiv u(\cdot,\omega)\quad\hbox{on $\R^N$} \qquad\hbox{for every $t>0$.}
\end{equation}
Take $\omega\in\Omega$ such that $w(\cdot,\omega)$ is a critical subsolution, $\widetilde C(\omega)\not=\emptyset$ 
and \eqref{Au1} holds. We  know that
$w(\cdot,\omega)\leq u(\cdot,\omega)$  in $\R^N$ and
$w(\cdot,\omega)= u(\cdot,\omega)$ in $\widetilde C(\omega)$.  
By the monotone character of the operator $\Tm$  we
deduce
\[
w(x,\omega)\leq (\Tm w)(x,\omega)\leq (\Tm u)(x,\omega)=w(x,\omega)\qquad\hbox{for every $t>0$.}
\]
Hence  all the inequalities in the
above formula  are indeed equalities, yielding 
$\widetilde C(\omega)\subseteq C(\omega)$ a.s. in $\omega$. 

In particular, we derive that 
$\A(\omega)=C(\omega)$ when $C(\omega)$ is almost surely empty, since the above argument implies, in this instance, that  stationary and nonempty closed random sets satisfying {\em (i)} in Definition \ref{def Aubry} do not exist. 

Let us then assume that $C(\omega)\not=\emptyset$ a.s. in $\omega$ and let us check that it satisfies  {\em (i)} in Definition \ref{def Aubry}. Pick a critical subsolution $g\in S$ and let $u$ be the admissible critical subsolution defined through \eqref{def u}. To prove that $u(\cdot,\omega)$ is an almost sure critical solution on $\R^N$, we only have to check, in view of \eqref{eq maximal subsol}, that the supersolution test is satisfied on $C(\omega)$ a.s. in $\omega$. But this follows in view of \eqref{eq prop2 w} and of Propositions \ref{prop1 w} and \ref{prop supersol on Xv} with $v:=u$. 
\qed

\medskip

\begin{oss}\label{oss A check}
Since  $v$ is a subsolution of $H=a$ if and only if $-v$ is a
subsolution of $\check H=a$, then  $H$ and $\check H$ have the
same critical value. Moreover, if $v$ is a critical subsolution
for $H$ which is weakly strict outside some stationary closed
random set $C(\omega)$, then $-v$ is critical for $\check H$ and
weakly strict outside $C(\omega)$, see Proposition 5.10 in
\cite{DS12} for more detail. In view of Theorem \ref{teo Aubry
set} and Remark \ref{oss minimal wstrict}, this implies that $H$
and $\check H$ have the same Aubry set.
\end{oss}

\bigskip

We  finally record for later use that, as a consequence of Theorem
\ref{teo Aubry set}, we are able  to extend Theorem 5.9 in
\cite{DS12}, employing the same argument used there, as follows:

\begin{teorema}\label{criti}  Assume that $\A(\omega)\not=\emptyset$ a.s. in $\omega$.
Then there exists a set $\widehat\Omega$ of probability 1 such
that for any
   $\omega\in\widehat\Omega$ and any $x\in\A(\omega)$ we can find
a curve $\eta_x:\R\to\A(\omega)$ (depending on $\omega$) with
$\eta_x(0)=x$ satisfying the following properties:
\begin{itemize}
\item[\em (i)] \quad for every $a<b$ in $\R$
\[
S(\eta_x(a),\eta_x(b),\omega)
    =
    \int_a^b L(\eta_x,\dot\eta_x,\omega) \,\dd s;
\]
\item[\em (ii)] \quad for every $v\in\S$ there exists a set
$\Omega_v$ of probability 1 such that for every
$\omega\in\Omega_v$
\[
 \int_a^b  L(\eta_x,\dot\eta_x,\omega) \,\dd s
    =
    v(\eta_x(b),\omega)-v(\eta_x(a),\omega)\qquad\hbox{for every $a<b$ in $\R$.}\medskip
\]
\end{itemize}
If condition \eqref{C} holds, we furthermore have  \quad
$\displaystyle{\lim_{t \rightarrow \pm \infty} |\eta_x(t)|= +
\infty.}$
\end{teorema}

\medskip
\begin{oss}\label{oss minimal wstrict} From Theorems \ref{teo Aubry set} and \ref{criti} we deduce that 
$\A(\omega)$ is the
minimal stationary closed random set for which there exists a
critical admissible subsolution which is weakly strict in its
complement.
\end{oss}

\medskip

\medskip

\end{subsection}

\begin{subsection}{Strict critical subsolutions.}   The purpose of this section is  to reinforce Theorem \ref{teo Aubry set}
showing the existence of a critical subsolution enjoying  the
property of being strict  outside the Aubry set in a
 stronger and more classical sense.
\begin{definition}
Let $C(\omega)$ a stationary closed random set. A critical subsolution $v\in\S$ is said to be {\em strict} in
$\R^N\setminus C(\omega)$ if the following property holds a.s. in $\omega$:\smallskip\\
\indent for every open set $U$ compactly contained in
$\R^N\setminus C(\omega)$ there exists $\delta>0$ such that
\begin{equation}\label{eq strict}
H(x,Dv(x,\omega),\omega)\leq -\delta\qquad\hbox{for a.e. $x\in
U$.}\medskip
\end{equation}

We will say that $v$ is { (weakly) strict}, with no further specification, to mean that it is (weakly) strict in $\R^N\setminus\A(\omega)$.
\end{definition}

\medskip

Next lemma makes precise that the previous  notion is actually a
strengthening of that of weakly strict subsolution. It is a purely
 deterministic result, where $\omega$ plays just the role of a  parameter, and so is  omitted for notational simplicity.

\smallskip

\begin{lemma}\label{lemma wstrict}
Let $v$ be a critical subsolution satisfying \eqref{eq strict}
in some open subset $U$ of $\R^N$ and for some $\delta>0$. Then
\begin{equation}\label{claim wstrict}
v(x)-v(y)<S(y,x)\qquad\hbox{for every $x\in U$ and $y\not=x$.}
\end{equation}
\end{lemma}

\begin{dimo}
Let $x\in U$ and $B_r(x)$ be any closed ball contained in $U$.
Since $S$ is a geodesic--type semidistance, it will be enough to check \eqref{claim wstrict} for every $y\in\partial B_r(x)$.

Let $\rho=\rho(0,-\delta)>0$ be chosen according to \eqref{spessore} and choose a curve
$\gamma$ joining $y$ to $x$ such that
\begin{equation}\label{eq quasi minimal}
\frac{\rho\,r}{2}+S(y,x)>\int_0^1
\sigma(\gamma,\dot\gamma)\,\dd s.
\end{equation}
Let us set $\tau:=\sup\{t\in [0,1]\,:\,\gamma(t)\in\R^N\setminus
B_r(x)\,\}$ and $z=\gamma(\tau)$. By taking into account
\eqref{spessore} and the fact that the critical subsolution
$v$  satisfies \eqref{eq strict} in $B_r(x)$, we
get:
\begin{eqnarray*}
    \int_0^1 \sigma(\gamma,\dot\gamma)\,\dd s
    &=&
    \int_0^\tau \sigma(\gamma,\dot\gamma)\,\dd s
    +
    \int_\tau^1 \sigma(\gamma,\dot\gamma)\,\dd s\\
    &\geq&
    v(z)-v(y)+v(x)-v(z)+\rho\,\int_\tau^1 |\dot\gamma(s)|\,\dd s.
\end{eqnarray*}
hence by \eqref{eq quasi minimal}
\[
S(y,x)\geq
v(x)-v(y)+\frac{\rho\,r}{2}>v(x)-v(y).
\]
\end{dimo}

\smallskip

The converse implication does not hold. More generally, the inequality
\[
 v(x)-v(y)<S_a(y,x)\qquad\hbox{for every $x,\,y$ in an open set $U\subseteq\R^N$}
\]
does not imply ${\D{ess}\sup}_{x\in U}\, H(x,Dv(x))<a$. For
instance, the antiderivative of a function vanishing on an open
and dense subset of $\R$ of small measure and  equal to 1 in the
complement is a weakly strict subsolution of the 1--dimensional
Eikonal equation $|u'|=1$ in $\R$, but it is not strict in
$\R$.\bigskip

\medskip

 The  statement of the main theorem reads as   follows:

\begin{teorema}\label{teo strict subsol}
Let $H$ satisfy (H1)--(H4). Then there exists a strict critical subsolution in $\S$.

More precisely, for every weakly strict critical subsolution
$w\in\S$ and every $\eps>0$, there exists a strict critical
subsolution $w_\eps\in\S$ such that
\begin{itemize}
 \item[\em (i)]\quad $\|w_\eps(\cdot,\omega)-w(\cdot,\omega)\|_\infty <
\eps\qquad\hbox{for every $\omega\in\Omega$;}$\smallskip
  \item[\em (ii)]\quad\  $w_\eps(\cdot,\omega)=w(\cdot,\omega)\quad\qquad\qquad\hbox{on $\A(\omega)$ \quad a.s. in $\omega$.}$
\end{itemize}
\end{teorema}

\smallskip

It, in particular, implies the existence of a critical admissible
subsolution, strict on the whole $\R^N$, when the random Aubry set
is almost surely empty.

\smallskip

As a consequence, we  also get:

\begin{cor}\label{cor strict subsol}
Let $H$ satisfy (H1)--(H4). The set of admissible, strict critical subsolutions is dense in $\S$ with respect to the Ky Fan metric on $L^0(\Omega;\D C(\R^N))$.
\end{cor}

\begin{dimo}
According to Theorems \ref{teo Aubry set} and \ref{teo strict subsol}, there exists a strict critical subsolution in $\S$, say it $v$. Now pick $u\in\S$ and set
\[
v_n(x,\oo)=\frac{1}{n}v(x,\omega)+(1-\frac{1}{n})u(x,\omega)\qquad{(x,\omega)\in\R^N\times\Omega}
\]
for every $n\in\N$. By convexity of the Hamiltonian, $v_n$ are strict critical subsolutions belonging to $\S$. Moreover
\[
v_n(\cdot,\omega)\ucv u(\cdot,\omega)\qquad\hbox{for every $\omega\in\Omega$,}
\]
meaning that $d(v_n(\cdot,\omega),u(\cdot,\omega))\to 0$ for every $\omega\in\Omega$. Since almost sure convergence implies convergence in probability, we get that $v_n$ converge to $u$ with respect to the Ky Fan metric in $L^0(\Omega;\D C(\R^N))$ in view of Theorem \ref{teo Ky Fan}.
\end{dimo}
\vspace{2ex}\\
\indent To pass from the existence of  a weakly strict admissible
subsolution to that of a strict one,
 we make use of two, in a sense complementary, crucial properties of the Lax--Oleinik semigroups that will be proved below.
 The first  is   the invariance of the values of  any critical subsolution  on
 the random Aubry set   under the action of $\Tm$ and $\Tp$,  the
 second  instead the strict monotonicity of  $\Tm$, when applied to a weakly strict critical subsolution,
 outside such set, at least for small times.

\smallskip

\begin{prop} \label{prop aubryconst}
Let $w \in \S$. Then the following property holds a.s. in $\omega$:
\[\Tm w(x,\oo)= \Tp w(x,\oo) =w(x,\oo) \quad\text{ for
any $x \in \A(\oo)$ and $t >0$.}\]
\end{prop}

\begin{dimo} We assume that Aubry set is a.s. nonempty otherwise
the statement is void. We take $\oo$ such that $\A(\oo) \neq
\emptyset$, $w(\cdot,\oo)$ is a subsolution of the corresponding
critical Hamilton--Jacobi equation and assertion {\em (ii)} of Theorem \ref{criti} holds for $v:=w$. Pick a point $x\in\A(\omega)$ and a time $t>0$.
According to Proposition \ref{prop critical sol},
\[
 \Tm w(x,\omega)\geq w(x,\omega).
\]
Now, let $\eta_x:\R\to\R^N$ be the curve chosen according to Theorem \ref{criti}. Then
\[
  w(x,\omega)
  =
  w(\eta_x(-t),\omega)+\int_{-t}^0  L(\eta_x,\dot\eta_x,\omega) \,\dd s
  \geq
  (\Tm w)(x,\omega),
\]
yielding equality. The assertion for $\Tp w$ can be proved analogously in view of Remarks \ref{oss H check} and \ref{oss A check}.
\end{dimo}

\medskip

\bigskip

\begin{prop}\label{lemma increasing}
Let $w\in\S$ be weakly strict. Then we can determine a set $\Omega_w$ of probability 1 such
that  for every $\oo \in \Oo_w$ and
$x\in\R^N\setminus\A(\omega)$ there exists $t_x=t_x(\omega)\in (0,+\infty]$ such that the
function
\[
t\mapsto v(t,x,\omega):=\left(\Tm w\right)(x,\omega)
\]
is strictly increasing in $[0,t_x)$ and constant for $t\geq t_x$.
\end{prop}

\begin{dimo}
Let $\Omega_w$ be a set of probability 1 made up by elements $\omega$ for which $w(\cdot,\oo)$ is a weakly strict critical subsolution. Fix $\omega\in\Omega_w$ and $x \in \R^N \setminus \A(\oo)$. Let us define $J$ as the set of times $t>0$ satisfying the following property:
 \begin{equation}\label{incr1}
(\Tm w)(x,\oo) =  w(y, \oo) + h_t(y,x,\oo)\qquad\hbox{for some $y \in \A(\oo)$}.
\end{equation}
Let us set $t_x:=\inf J$, where we agree that $t_x=+\infty$ if $J$ is empty.
From Proposition  \ref{prop B4}, which bounds the distance of the point $y$ in \eqref{incr1} from $x$, and the fact that $\A(\omega)$ is closed we infer that $J$ is closed and that $t_x>0$.

Let us prove that $t\mapsto \left(\Tm w\right)(x,\omega)$ is
strictly increasing in $[0,t_x)$.  Take $s >t$ in $[0,t_x)$, then
\begin{equation}\label{incr2}
   (T^-_s\, w)(x,\omega)=  (T^-_{s - t}\, w)(y,\omega) + h_t(y,x,\oo) \qquad
\hbox{for some  $y \in \R^N$.}
\end{equation}
If $y \in \A(\oo)$, from Proposition \ref{prop aubryconst} we deduce
\[
(T^-_s\, w)(x,\omega)=   w(y,\omega) + h_t(y,x,\oo) > (\Tm
w)(x, \oo),
\]
were the strict inequality comes from the fact that $t\not\in J$.
If instead  $y \not \in
\A(\oo)$,  we invoke Lemma \ref{lemma pre increasing}, which holds
true for the $\oo$ we are working with, to get from \eqref{incr2}
\[
(T^-_s\, w)(x,\omega) >   w (y,\omega) + h_t(y,x,\oo) \geq (\Tm
w)(x,\oo).
\]

Let us now prove that $ t \mapsto (\Tm w)(x,\oo)$ is constant in $[t_x, + \infty)$ when $t_x<+\infty$.
Let $y\in\A(\omega)$ be a point satisfying \eqref{incr1} with $t_x$ in place of $t$.
We invoke Proposition \ref{prop aubryconst} to get  for  $t >t_x$
\[
(\Tm w)(x,\oo) \leq  (T^-_{t-t_x} w)(y, \oo) + h_{t_x}(y,x,\oo) = w(y, \oo) +
h_{t_x}(y,x,\oo)= (T^-_{t_x} w)(x,\oo).
\]
By  monotonicity
properties of Lax-Oleinik semigroup pointed out in Proposition
\ref{prop critical sol}, we get the assertion.
\end{dimo}

\bigskip

To proceed in our analysis, we  need  two technical lemmata about
locally Lipschitz functions and related Clarke's generalized
gradients. In what follows, we denote by $\pi_1$, $\pi_2$ the maps
defined as $\pi_1(p_t,p_x)=p_t$, $\pi_2(p_t,p_x)=p_x$ for every
$(p_t,p_x)\in\R\times\R^N$.

\smallskip

\begin{lemma}\label{lemma tecnico1}
Let $v(t,x)$ be a locally Lipschitz function in $(0,+\infty)\times\R^N$ such that, for every bounded open subset $U$
of $\R^N$, the functions
\[
\{v(\cdot,x)\,:\,x\in U\,\}\quad\hbox{are locally equi--semiconcave (resp. equi--semiconvex) in $(0,+\infty)$.}
\]
Then, for every $(t_0,x_0)\in (0,+\infty)\times\R^N$,
\[
 \pi_1\left(\partial^c v(t_0,x_0)\right)=\partial^c_t v(t_0,x_0),\qquad \pi_2\left(\partial^c v(t_0,x_0)\right)\supseteq\partial^c_x v(t_0,x_0).
\]
In particular, the map $(t,x)\mapsto \partial^c_t v(t,x)$ is upper semicontinuous in $(0,+\infty)\times\R^N$.
\end{lemma}

\begin{dimo}
The fact that
\[
\pi_1\left(\partial^c v(t_0,x_0)\right)\supseteq\partial^c_t v(t_0,x_0),\qquad \pi_2\left(\partial^c v(t_0,x_0)\right)\supseteq\partial^c_x v(t_0,x_0)
\]
for every $(t_0,x_0)\in (0,+\infty)\times\R^N$ follows from Proposition 2.3.16 of \cite{Cl}. To prove that
$\pi_1\left(\partial^c v(t_0,x_0)\right)=\partial^c_t v(t_0,x_0)$, it will be enough to show that
\[
\pi_1\left(\partial^* v(t_0,x_0)\right)\subseteq\partial^c_t v(t_0,x_0).
\]
Let $p_{t_0}\in \pi_1\left(\partial^* v(t_0,x_0)\right)$ and let $(t_n,x_n)$ be a sequence of differentiability points for $v$ converging to $(t_0,x_0)$ such that $\partial_t v(t_n,x_n)=:p_{t_n}\to p_{t_0}$ as $n\to +\infty$. The functions
\[
 \phi_n(t):=w(t-t_0+t_n,x_n)
\]
are locally equi--semiconcave (resp. equi--semiconvex) in $(0,+\infty)$  and
\[
\phi_n \ucv v(\cdot,x_0) \quad \hbox{in $(0,+\infty)$}.
\]
Moreover $\phi_n$ are differentiable at $t_0$ and $\phi'_n(t_0)=p_{t_n}$, in particular $p_{t_n}$ belongs to the superdifferential (resp. subdifferential) of $\phi_n$ at $t_0$ by semiconcavity (resp. semiconvexity). This property passes to the limit since the functions $\phi_n$ are locally equi--semiconcave (resp. equi--semiconvex), hence
$p_{t_0}$ belongs to the superdifferential (resp. subdifferential) of the function $v(\cdot,x_0)$ at $t_0$, that is to $\partial^c_t v(t_0,x_0)$.

Since the map $(t,x)\mapsto \partial^c v(t,x)$ is upper semicontinuous in $(0,+\infty)\times\R^N$, so is the map
$(t,x)\mapsto \pi_1(\partial^c v(t,x))$, from which we get the last assertion in view of the identity just established.
\end{dimo}

\medskip

\begin{lemma}\label{lemma tecnico2}
Let $v(t,x)$ be a locally Lipschitz function in $(0,+\infty)\times\R^N$ and assume there is a bounded open subset $U$
of $\R^N$ such that
\begin{itemize}
 \item[\em (i)] the functions $\{v(\cdot,x)\,:\,x\in U\,\}$ are locally equi--semiconcave (resp. equi--semiconvex) in $(0,+\infty)$;\smallskip
  \item[\em (ii)] for every $x\in U$ the map
\[
 t\mapsto v(t,x)\qquad\hbox{is strictly increasing in $[0,t_x)$ }
\]
for some $t_x\in (0,+\infty]$.
\end{itemize}
Then, for every $x\in U$, the set
\[
 I_x:=\{t>0\,:\,\min\{p_t\,:\,p_t\in\partial^c_t v(t,x)\}>0\,\}
\]
is open and dense in $(0,t_x)$.
\end{lemma}

\begin{dimo}
Let us fix $x \in U$ and set
\[
f(t) = \min \{r \, : \, p_t \in \p_t v(t,x)\}\qquad\hbox{for every $t>0$.}
\]
Thanks to Lemma \ref{lemma tecnico1}, the set--valued map $t\mapsto\p_t v(t,x)$ is
upper semicontinuous. Consequently, the function $f$ is lower semicontinuous, so
that
\[
I_x= \{ t >0 \, : \, f(t) >0\}
\]
is open. Let us prove that it is dense in $(0,t_x)$. This
is indeed a consequence of the strict monotonicity of $t \mapsto
v(t,x)$ in $[0,t_x)$. Because of it, in fact, in any
subinterval of $(0, t_x)$ we find differentiability points of
$v(\cdot,x)$ with positive derivative. Since any such
point, say $t_0$, is of strict differentiability due to the
semiconcavity (resp. semiconvexity) of $v(\cdot,x)$, the corresponding derivative is also the unique
generalized gradient, which shows that $t_0 \in I_x$ and proves
the statement.
\end{dimo}

\bigskip

Our strategy for proving  Theorem \ref {teo strict subsol} is
structured in two parts: we first show the result under the
additional assumption that the $H$ is strictly convex and then
generalize it to Hamiltonians solely convex by means of a
regularization in time of the action of the negative Lax--Oleinik
semigroup via t--partial sup convolutions.

\bigskip

\noindent{\bf Proof of Theorem \ref{teo strict subsol} for $H$
strictly convex.}  The precise statement of our extra assumption
is:
\begin{itemize}
 \item[(H2$'$)] \quad$H(x,\cdot,\omega)\ \hbox{is strictly convex on $\R^N$}$\ \ \qquad for every $(x,\omega)\in\R^N\times\Omega$.
\end{itemize}
Let $w\in\S$ be a weakly strict subsolution and let $\eps>0$ be
fixed.  The main effect of   (H2$'$) is that the function
\begin{equation}\label{def v(t,x)}
 v(t,x,\omega):=(\Tm w)(x,\omega),\qquad{(t,x,\omega)\in (0,+\infty)\times\R^N\times\Omega,}
\end{equation}
is locally semiconcave in $(0,+\infty)$ with respect to $t$, see Lemma 2.11 in \cite{DavZav}. More precisely, for every open and bounded set $U\subset\R^N$ and every fixed $\omega\in\Omega$ the functions
\[
 \{v(\cdot,x,\omega)\,:\,{x\in U}\}\quad\hbox{are locally equi--semiconcave in $(0,+\infty)$.}
\]
Let $\kappa$ be the constant given by \eqref{def kappa_a} with $a=0$ and let $R(\kappa)$ be chosen according to Proposition \ref{prop B4}. Choose $\tau>0$ such that $\tau\, R(\kappa)<\eps$ and let $(t_n)_n$ be a dense sequence in $(0,\tau)$. We define
\begin{equation}\label{def w-eps}
 w_\eps(x,\omega)=\sum_{n=1}^{+\infty} \frac{1}{2^n} \, v(t_n,x,\oo)
,\qquad{(x,\omega)\in\R^N\times\Omega.}
\end{equation}
By Proposition \ref{prop B subsol}, we get
\[
 \|w_\eps(\cdot,\omega)-w(\cdot,\omega)\|_\infty
  \leq
  \sum_{n=1}^{+\infty} \frac{1}{2^n} \,\|v(t_n,\cdot,\omega)-w(\cdot,\omega)\|_\infty
  \leq \tau R(\kappa) <\eps
\]
for every $\omega\in\Omega$, showing in particular that $w_\eps(\cdot,\omega)$ is finite--valued.
As a convex combination of admissible critical subsolutions, see Proposition \ref{prop B subsol},
a standard argument shows that $w_\eps$ is an admissible subsolution as
well. From Proposition \ref{prop aubryconst} we also infer
\begin{equation}\label{equality on Aubry}
 w_\eps(\cdot,\omega)=w(\cdot,\omega)\qquad\hbox{on $\A(\omega)$}
\end{equation}
a.s. in $\omega$.
It is left to show that $w_\eps$ is strict.  To this purpose, let us fix $\omega$ in a set of probability 1 made up of  elements for which the assertion of Proposition \ref{lemma increasing} holds true and the functions $v(t_n,\cdot,\omega)$ are critical subsolutions. Pick a point $y\in\R^N\setminus\A(\omega)$. According to Lemma \ref{lemma tecnico2} and by density of the sequence $(t_n)_n$ in $(0,\tau)$, there exist $i\in\N$ and $a>0$ such that
\begin{equation}\label{ineq local strict}
\min\{p_t\,:\,p_t\in\partial^c_t v(t_i,y,\omega)\}>a.
\end{equation}
By upper semicontinuity of the set--valued map $x\mapsto \partial^c_t v(t_i,x,\omega)$, we infer that \eqref{ineq local strict} keep holding in a ball $B_r(y)$. Now we exploit the fact that $v(\cdot,\cdot,\omega)$ is a (sub)solution of the time--dependent equation \eqref{Cauchy problem}, in view of Proposition \ref{prop Tm solution}, so
\[
 p_t+H(x,p_x,\omega)\leq 0\qquad\hbox{for every $(p_t,p_x)\in\partial^c v(t_i,x,\omega)$ and $x\in B_r(y)$.}
\]
In view of Lemma \ref{lemma tecnico1} and of what remarked above, we get in particular
\[
 H(x,D_x v(t_i,x,\omega),\omega)<-a<0\qquad\hbox{for a.e. $x\in B_r(y)$.}
\]
By the definition of $w_\eps$ and the convexity of $H$, we conclude that
\begin{eqnarray*}
H(x,Dw_\eps(x,\omega),\omega)
    &\leq&
    \frac{1}{2^i}\,H (x,D_x v(t_i,x,\omega),\oo)\\
    &+&
    \sum_{n\not=i} \frac{1}{2^n}\,H (x,D_x v(t_i,x,\omega),\oo)
    \leq - \frac{a}{2^i}<0
\end{eqnarray*}
for a.e. $x\in B_r(y)$. This  actually shows that $w_\eps$ is strict since $y$ was arbitrarily
chosen in $\R^N\setminus\A(\omega)$.
\qed

\medskip

Looking carefully at the above argument, we recognize that
definition \eqref{def v(t,x)} can be interpreted as a convenient
way to select a 1--parameter family $\{v(t,\cdot,\cdot)\}_{t> 0}$
 of elements in $\S$ in such a way that the following conditions are satisfied almost surely: $v(\cdot,\cdot,\omega)$
 is a subsolution of the time--dependent equation \eqref{Cauchy problem};  the function $t\mapsto v(t,x,\omega)$
 is constant on $\A(\omega)$, and strictly increasing outside $\A(\omega)$ in a suitable neighborhood $[0,t_x)$ of $t=0$;
 the map $t\mapsto v(t,x,\omega)$  is locally semiconcave (or semiconvex) in $(0,+\infty)$, with a modulus that is locally uniform with respect to $x$.

In the general case of  an Hamiltonian, not strictly convex, but
just convex, the latter condition is apparently no longer
fulfilled by the random variable given by \eqref{def v(t,x)}.To get the proof of Theorem \ref{teo strict subsol} in
full generality, we   modify the definition of $v$ by setting
\begin{equation}\label{convo}
    v(t,x,\oo)= \sup_{s \geq 0} \left \{ (T^-_s w)(x,\oo) -
\frac1{2\,\de}\, (s -t)^2 \right \}.
\end{equation}
To complete our task it is then enough to check  out  that   such
a $v$ fulfills {\em all} the requirements listed above,  and this
is indeed the content of our next proposition. In this way $v$ can
be actually used in formula \eqref{def w-eps} to provide a strict
subsolution $w_\eps\in\S$ almost surely satisfying \eqref{equality
on Aubry}, while for the inequality $\|w_\eps-w\|_\infty<\eps$ in
the item {\em (i)} of the statement, it  simply suffices to choose
$\tau>0$ and $\delta>0$ small enough.

\medskip

\begin{prop}
Let $v$ be defined via \eqref{convo} for some $w\in\S$ and $\delta>0$. The following properties hold:\smallskip
\begin{itemize}
 \item[\em (i)] $\|v(t,\cdot,\omega)-w(\cdot,\omega)\|_\infty\leq R(\kappa)(t+2\delta R(\kappa))$\quad for every $t>0$ and $\omega\in\Omega$;\smallskip
 \item[\em (ii)] the function $v(\cdot,\cdot,\omega)$ is a Lipschitz subsolution in $(0,+\infty)\times\R^N$ of the time--dependent equation \eqref{Cauchy problem} for every $\omega\in\Omega$;\smallskip
 \item[\em (iii)] $v(t,\cdot,\cdot)\in\S$\quad for every $t>0$;\smallskip
 \item[\em (iv)] for every $\omega$ in a set of probability 1 the following property holds:
\[
 v(t,\cdot,\omega)=w(\cdot,\omega)\quad\hbox{on $\A(\omega)$\qquad for every $t>0$;}
\]
\item[\em (v)] the functions $\{v(\cdot,x,\omega)\,:\,x\in\R^N,\ \omega\in\Omega\,\}$ are equi--semiconvex in $(0,+\infty)$.\medskip
\end{itemize}
Moreover, if $w\in\S$ is weakly strict, then $v$ enjoys the assertion of Proposition \ref{lemma increasing}.
\end{prop}

\begin{dimo}
By Proposition \ref{prop B4} we know that the function
$(t,x) \mapsto (\Tm w)(x,\oo)$ is $R(\kappa)$--Lipschitz continuous
$[0,+\infty)\times \R^N$ for every $\oo\in\Omega$, in particular $(\Tm
w)(x,\oo)$ grows at most linearly in $t$. This implies that the supremum in \eqref{convo} is attained. Let us denote by $Y(t,x,\omega)$ the set of maximizers of the right--hand side term of \eqref{convo}.
Then
\[
 |s-t|\leq 2\delta R(\kappa)\qquad\hbox{for every $s\in Y(t,x,\omega)$.}
\]
This follows from the fact that, for any such $s$, we have
\[
  \frac 1{2 \, \de} (t-s)^2 \leq  (T^-_s w)(x,\oo) - (\Tm w)(x,\oo)\leq R(\kappa)|t-s|.
\]
In particular,
\[
0\leq v(t,x,\omega)-(\Tm w)(x,\omega)\leq (T^-_s w)(x,\oo) - (\Tm w)(x,\oo)\leq 2\delta R(\kappa)^2,
\]
and assertion {\em (i)} follows by Proposition \ref{prop B4}--{\em (ii)}.

Assertion {\em (ii)} is well known, see for instance \cite{CaSi}: the Lipschitz character of $v(\cdot,\cdot,\omega)$ comes from the fact that it is the supremum of equi--Lipschitz functions; the subsolution property is a consequence of
the inclusion
\[
D^+ v(t,x,\omega)\subseteq D^+ (T^-_{s}w)(x,\omega)\qquad\hbox{for every $s\in Y(t,x,\omega)$,}
\]
together with the fact that $(\Tm w)(x,\omega)$ is a (sub)solution of \eqref{Cauchy problem} and that the Hamiltonian is autonomous.

Let us prove item {\em (iii)} for any fixed $t>0$. First notice that, by continuity with respect to $s$, the supremum in \eqref{convo} can be taken over a dense and countable subset of $(0,+\infty)$, yielding that $v(t,\cdot,\cdot)$ is a random variable. Moreover,  $v(t,\cdot,\omega)$ is the supremum of a family of critical subsolutions whenever $w(\cdot,\omega)$ is a critical subsolution, see Proposition \ref{prop critical sol}, i.e. almost surely. This implies that $v(t,\cdot,\omega)$ is a critical subsolution a.s. in $\omega$. The fact that $v(t,\cdot,\omega)$ is almost surely sublinear is obvious from {\em (i)} and from the fact that $w\in\S$. It is left to show that $v(t,\cdot,\cdot)$ has stationary increments. Indeed,
$w$ has stationary increments, meaning that, for every fixed $z\in\R^N$, there exists a random variable $k(\oo)$ and a
set $\Oo_z$ with probability $1$ such that
\[
w(\cdot + z, \oo) = w(\cdot, \tau_z \oo) + k(\oo) \quad\hbox{on $\R^N$}\qquad\hbox{for every $\omega\in\Omega_z$.}
\]
Keeping in mind the stationary character of the Lagrangian, it is easily checked that the same identity is satisfied by
$T^-_s w$ for every $s>0$, and, consequently, by $v(t,\cdot,\cdot)$.

Assertion {\em (iv)} is straightforward consequence of the definition of $v$ in view of Proposition \ref{prop aubryconst}.

Assertion {\em (v)} is also well known, see for instance \cite{CaSi}.\\
Let us prove the last assertion.
We first note that, when $w(\cdot,\omega)$ is a critical subsolution, the monotonicity of the map $t\mapsto(\Tm w)(x,\omega)$ readily implies, by the very definition of $v$, that
\begin{equation}\label{eq Y}
Y(t,x,\omega)\subset [t,+\infty).
\end{equation}
Assume now that $w$ is weakly strict and let us denote by $\Omega_w$ a set of probability 1 made up by elements $\omega$ for which $w(\cdot,\oo)$ is a weakly strict critical subsolution. Fix $\omega\in\Omega_w$ and pick a point $x\in\R^N\setminus\A(\omega)$. Then we know by (the proof of) Proposition \ref{lemma increasing} that there exists $t_x\in (0,+\infty]$ such that the function $t\mapsto(\Tm w)(x,\omega)$ is strictly increasing in $[0,t_x)$ and constant for $t\geq t_x$. In view of \eqref{eq Y} and the very definition of $v$, we infer that
\[
 t\mapsto v(t,x,\omega)\qquad\hbox{is constant for $t\geq t_x$.}
\]
It is left to show that it is strictly increasing in $[0,t_x)$.
Take $t_1> t_2$  in $[0,t_x)$ and let $s_i\in Y(t_i,x,\omega)$ for $i=1,\,2$. Two cases are possible: either
$s_2< t_x$ or $s_2=t_x$. In the first instance we have
\begin{eqnarray*}
  v(t_2,x ,\oo) &=& (T^-_{s_2}w)(x,\oo) - \frac1{2\,\de}\, (s_2-t_2)^2\\
&<&\big(T^-_{t_1+(s_2 -t_2)} w\big)(x,\oo) - \frac1{2\,\de}\, (s_2-t_2)^2
  \leq v(t_1,x ,\oo),
\end{eqnarray*}
in the second
\begin{eqnarray*}
  v(t_2,x ,\oo) &=& (T^-_{t_x}w)(x,\oo) - \frac1{2\,\de}\, (t_x - t_2)^2 \\
   &<&  (T^-_{t_x}w)(x,\oo)  - \frac1{2\,\de}\, (t_x - t_1)^2 \leq
   v(t_1,x ,\oo),
\end{eqnarray*}
as it was to be shown.
\end{dimo}

\end{subsection}
\end{section}

\begin{section}{Tonelli Hamiltonians}\label{Tonelli}

In this section we deal with Hamiltonians  satisfying more
stringent regularity assumptions and named after Tonelli. In the
first subsection we provide basic definitions and illustrate the
salient features of the corresponding  Hamiltonian flow and
Lax--Oleinik semigroups. In the second one we prove  existence of
$C^{1,1}$ strict subsolutions in the stationary ergodic setting
and investigate their properties. This is achieved by applying
Bernard's method on $\CC^{1,1}$--regularization of strict
subsolutions on compact manifolds in the deterministic case.
\smallskip

Throughout the section we will use the term  semiconcave
(respectively semiconvex) in a stronger sense than the one defined
in Section \ref{sez basic}:  we will  in fact additionally require 
the modulus  to be linear, namely  $\Theta(h)=K\,h$ for
some $K >0$. If such a constant need to be showcased  then  we
will employ the diction $K$--semiconcave (respectively,
$K$--semiconvex).   We recall that a function $u$ is both
$K$--semiconcave and $K$--semiconvex in  some open subset $U$ of
$\R^N$ if and only if is of class $C^{1,1}$ in $U$ and
$\D{Lip}(Du;U)\leq K$, see \cite{CaSi}.

\begin{subsection}{Generalities.}\label{Tonelligen}

We say that a stationary ergodic Hamiltonian $H$ is {\em Tonelli}
if it enjoys conditions (H1)--(H4) and the following set of assumptions:
\begin{itemize}
    \item[(T1)] \quad $H(\cdot,\cdot,\omega)\in C^2(\R^N\times\R^N)\quad$ for every $\omega\in\Omega$;\medskip
    \item[(T2)] \quad for every $R>0$ there exists a constant $\nu_R>0$ such that
\[
    \frac{\partial^2 H}{\partial p^2}(x,p,\omega)> \nu_R\,\D{Id}_{\R^N}\quad\hbox{for every
    $(x,p,\omega)\in\R^N\times B_R\times\Omega$;} \medskip
\]
     \item[(T3)] \quad for every $R>0$ there exists $L_R>0$ such that
     \begin{eqnarray*}
     \|DH(\cdot,\cdot,\omega)\|_{L^\infty(\R^N\times B_R)}&\leq&
     L_R,\medskip\\
     \D{Lip}\left(DH(\cdot,\cdot,\omega);\,\R^N\times
     B_R\right)&\leq&
     L_R\vspace{2ex}
     \end{eqnarray*}
     \quad for every $\omega\in\Omega$.\\
\end{itemize}

Under above assumptions, it is well known that the associated
Lagrangian is of class $C^2$ as well, see for instance
\cite{CoIt00,Fathi}, and satisfies properties analogous to
(T1)--(T3).

\begin{oss}\label{tonton} The above hypotheses, in particular (T2) and (T3),  are  adaptation to the
stationary ergodic environment of the usual ones required for
deterministic Tonelli Hamiltonians. The changes are basically
due to the fact that we need bounds independent of $\oo$,  and this
immediately implies that they have also to be global in $x$. In
fact, bounds independent of $\omega$ that are local in $x$ simply do not make sense in our
frame, for stationarity and ergodicity assumptions should
automatically transfer them to the whole $\R^N$. Similarly,  we
could rephrase (T2) and (T3)  by requiring $\nu_R$ and $L_R$ to
depend in a measurable way on $\omega$: the  ergodicity assumption
would then imply that they are almost surely constant.
\end{oss}

We will denote by $\phi^H_t(x,p,\omega)$ the Hamiltonian flow,
i.e. the flow associated with Hamilton's ODE
\begin{eqnarray*}
  \begin{cases}
    \dot \xi=\ \ \frac{\partial H}{\partial p}(\xi,\eta,\omega)\medskip\\
    \dot  \eta=-\frac{\partial H}{\partial x}(\xi,\eta,\omega).
  \end{cases}
\end{eqnarray*}
 The corresponding integral curves will be also called  {\em characteristics} in the sequel. As well known,
  $H(\phi^H_t(x,p,\omega),\omega)=H(x,p,\omega)$ for
every $(x,p,\omega)$, which yields, by the coercivity assumption
(H3), that the flow is complete, i.e. globally defined in time.

\smallskip

 We proceed discussing the two main additional features of Lax--Oleinik semigroups
 for Tonelli Hamiltonians, namely the fact that the action of
 semigroups is driven by characteristics and secondly that the
 functions obtained  in this way are semiconcave/semiconvex. See
\cite{Fathi} for corresponding proofs.

\medskip

\begin{prop}\label{prop minimizer}
Let $u$ be a $\vartheta$--Lipschitz  random function, $x\in\R^N$
and $t>0$. Let $\gamma:[0,t]\to\R^N$ be a Lipschitz curve with
$\gamma(t)=x$  such that
\[
(T_t^-u)(x,\omega)=u(\gamma(0),\omega)+\int_{0}^t
L(\gamma,\dot\gamma,\omega)\,\dd s.
\]
%
Then $\gamma$ is  the projection on the state variable space of a
characteristic taking value  $(\gamma(0), p_{\gamma(0)})$ at time 0, and
$(x,p_x)$ at time $t$, where
\[
p_{\gamma(0)}\in
  {D^-}u(\gamma(0),\omega),
\qquad
p_x\in {D^+}(T_t^-u)(x,\omega).
\]
\end{prop}

\medskip

While for continuous Hamiltonians, we can
only assert $\kappa$--Lipschitz continuity of $(T^-_t u)(\cdot,\oo)$ and
$(T^+_t u)(\cdot,\oo)$ when $u \in \S$, for Tonelli ones we get
semiconcavity and semiconvexity, respectively.  This
property will be crucial to transfer to the stationary ergodic
setting the regularization  procedure yielding $\CC^{1,1}$
critical subsolutions.\medskip

\begin{prop}\label{lemma ht}
Let $t_0>0$. Then the functions
$\{\,h_t(\cdot,\cdot,\omega)\,:\,\omega\in\Omega,\,t\geq t_0\,\}$
are locally equi--concave (and equi--Lipschitz) in
$\R^N\times\R^N$.
In particular, there exists a constant $K=K(t_0)$
such that, for every $w\in\S$ and for every $t\geq t_0$, \smallskip
\begin{eqnarray*}
(\Tm w)(\cdot,\omega)\ &\hbox{is}&\ \hbox{$K$--semiconcave in
$\R^N$,}\\
(\Tp w)(\cdot,\omega)\ &\hbox{is}&\ \hbox{$K$--semiconvex in
$\R^N$,}
\end{eqnarray*}
for every $\omega\in\Omega$.
\end{prop}

\begin{oss}\label{tonton2}
The fact that the constants of local semiconcavity and
semiconvexity are global in $\R^N$ and independent of $\oo$ 
is a consequence of the
fact that, for every $R
>0$, the Lipschitz constant of $DL(\cdot,\cdot,\oo)$ on $\R^N
\times B_R$ is finite and independent of $\oo$. 
Condition (T2) is crucial to get such a bound, as it can be deduced from the  known relation   
\[
\frac{\partial L}{\partial q}(x,q,\omega)= \left(\frac{\partial
H}{\partial p}\right)^{-1}(x,q,\omega).
\]
\end{oss}

\medskip
We can derive from Theorem \ref{criti}:

\smallskip

\begin{teorema}\label{teo differentiability}
The following facts hold:
\begin{itemize}
  \item[\em (i)] for every $v\in\S$, there exists a set
  $\Omega_v$ of probability 1 such that for every
  $\omega\in\Omega_v$ and every $x\in\A(\omega)$
  \[
    v(\cdot,\omega)\quad\hbox{is differentiable at $x$}\quad\hbox{and}\quad H(x,Dv(x))=0.
  \]
  \item[\em (ii)] Let $u$ and $v$ belong to $\S$. Then
  \[
    Du(\cdot,\omega)=Dv(\cdot,\omega)\qquad\hbox{on $\A(\omega)$}
  \]
  for every $\omega\in\Omega_u\cap\Omega_v$.
\end{itemize}
\end{teorema}

Theorem \ref{teo differentiability} is analogous to a well known
result in weak KAM Theory and can be proved similarly, see
\cite{Fathi}.

Let us pick $v\in\S$ and set
\[
\widetilde
\A(\omega):=\{(x,Dv(x,\omega))\,:\,x\in\A(\omega)\,\},\qquad\hbox{$\omega\in\Omega$.}
\]
Up to a set of null probability, the definition of $\widetilde
\A(\omega)$ is independent of the choice of $v\in\S$ in view of
Theorem \ref{teo differentiability}. The following holds, see
\cite{Fathi}:

\begin{teorema}\label{teo flow}
There exists a set $\widehat\Omega$ of probability 1 such that for
every $\omega\in\widehat\Omega$
\[
\phi_H^t\big(\widetilde\A(\omega),\omega\big)=
\widetilde\A(\omega)\qquad\hbox{for every $t\in\R$.}
\]
\end{teorema}
\medskip
\end{subsection}

\begin{subsection}{$\CC^{1,1}$ subsolutions in stationary ergodic
case.} \label{random}

In this subsection we show how to pass from a strict critical and admissible subsolution, which is in general just Lipschitz--continuous with respect to $x$, to one which is of class $C^{1,1}$ in $\R^N$.  The precise result we will prove is the following:

\begin{teorema}\label{sawbis}
Let $H$ satisfy (H1)--(H4) and
(T1)--(T3). Then there exists a strict critical subsolution in $\S$ of class $C^{1,1}$ in $\R^N$.

More precisely, for every strict critical subsolution
$w\in\S$ and every $\eps>0$, there exists a strict critical
subsolution $w_\eps\in\S$ such that
\begin{itemize}
 \item[\em (i)]\quad $\|w_\eps(\cdot,\omega)-w(\cdot,\omega)\|_\infty <
\eps\qquad\hbox{for every $\omega\in\Omega$;}$\smallskip
  \item[\em (ii)]\quad\  $w_\eps(\cdot,\omega)=w(\cdot,\omega)\quad\qquad\qquad\hbox{on $\A(\omega)$ \quad a.s. in $\omega$;}$\smallskip
  \item[\em (iii)]\quad\ $w_\eps(\cdot,\omega)\in\CC^{1,1}(\R^N)$ \qquad\quad\ \  {for every $\omega\in\Omega$}.
\end{itemize}
\end{teorema}

This theorem, in particular, implies the existence of a critical
admissible subsolution of class $\CC^{1,1}$ and strict on the whole $\R^N$ when the
random Aubry set is almost surely empty.

As a direct consequence of Theorem \ref{sawbis} and Corollary \ref{cor strict subsol} we get

\begin{cor}
Let $H$ satisfy (H1)--(H4) and (T1)--(T3). Then the set of admissible, strict critical subsolutions of class $C^{1,1}$ in $\R^N$ is dense in $\S$ with respect to the Ky Fan metric on $L^0(\Omega;\D C(\R^N))$.
\end{cor}

\medskip

We start by showing a further invariance property of Lax--Oleinik
semigroup.

\begin{prop}\label{invariance}
Let $H$ be a stationary ergodic Tonelli Hamiltonian.  Then the family of strict subsolution of $\S$ is invariant for the positive and negative Lax--Oleinik semigroups.
  Namely, if $u\in\S$ and is strict, then both $\Tp u$ and $\Tm u$
 belong to $\S$ and are strict, for every $t>0$.
\end{prop}
\begin{dimo}  Let $u \in \S$ be strict. Take $\oo$ such that
$u(\cdot,\oo)$ is a  strict  critical subsolution, $u(\cdot,\oo)$
is  differentiable on $\A(\oo)$ and $\widetilde \A(\oo)$ is
invariant with respect to the Hamiltonian flow. Theorems \ref{teo
differentiability} and \ref{teo flow} guarantee that such a choice
can be done in set of probability $1$.

 We go on fixing $t
>0$ and a closed ball $B$ disjoint from $\A(\oo)$. We consider
the set $V$ made up by points $y$ optimal for $(\Tm u)(x,\omega)$, for
some $x \in B$. It inherits from $B$ the property of being closed.
 As established in
Proposition \ref{prop minimizer},   if $y \in V$  is optimal for
$(\Tm u)(x,\omega)$, the corresponding Lagrangian minimizer $\gamma: [0,t]\to\R^N$
is the projection on the state variable space of the characteristic  taking
the value $(y,p)$ at time $0$, with $p \in D^-u(y,\oo)$.  If $y \in
\A(\oo)$, then $p= Du(y,\oo)$ and $(y,p) \in \widetilde
\A(\oo)$, but then, by the invariance property of $\widetilde \A(\oo)$,
$\gamma$ cannot reach $x$, which does not belong to $\A(\oo)$. This
argument shows by contradiction that $V \cap\A(\oo) =
\emptyset$.

Being $u(\cdot,\oo)$ strict critical subsolution, there thus
exists $\de >0$ with
\[H(y,  D u(y,\oo),\omega) \leq -\de \quad\text{ for a.e. $y \in V$.}
\]
Taking into account that the Hamiltonian stays constant on
characteristic, we derive
\[
H(x,D_x (\Tm u)(x,\oo),\omega) \leq -\de \quad\text{ for a.e. $x \in B$.}
\]
Since $B$ is a closed ball arbitrarily chosen in the complement of
$\A(\oo)$, this proves the assertion for $\Tm u$, and consequently for
$\Tp u$ in view of Remark \ref{oss H check}.
\end{dimo}

\medskip

The same argument of the above proposition yields a strengthened
form of Proposition  \ref{lemma increasing} for Tonelli
Hamiltonians. Namely:

\smallskip

\begin{cor}\label{cor increasing bis}
Let $H$ be a stationary ergodic Tonelli Hamiltonian. Then, for any  weakly
strict subsolution  $w\in\S$, there exists a  set $\Omega_w$ of probability 1 such
that,  for every $\oo \in \Oo_w$ and $x\in\R^N\setminus\A(\omega)$,  the function
\[
t\mapsto \left(\Tm w\right)(x,\omega)
\]
is strictly increasing in   $[0,+\infty)$. In particular, $\Tm w$ is weakly strict 
for every $t>0$.
\end{cor}

\medskip
For the sake of completeness,  we also  provide, in the framework
of the proof of Theorem \ref{sawbis}, a worked out presentation of
Bernard regularization technique in the deterministic setting
using our terminology and notations. The material is illustrated
in a more elementary way and providing  more details  with respect to
the original proof in \cite{Be}, at the expense of some loss of
concision and elegance.

The  first two steps are purely deterministic, and so  $\omega$ is
omitted. Lemma \ref{regular} establish that a function which is
$\CC^{1,1}$ locally around a given point, {\em remains} of class
$\CC^{1,1}$, at least in a smaller neighborhood of the same point,
under application  of $T^-_t$ for small times. It will be used in
the subsequent Proposition \ref{semi} to show, by working on
subtangents, that a function locally semiconvex {\em becomes}
locally $\CC^{1,1}$ under the same action.

\smallskip

In what follows, the symbols $\pi_1$ and $\pi_2$ will
denote the projections on the space  and momentum variable,
respectively, i.e. the maps defined as $\pi_1(x,p)=x$,  $\pi_2(x,p)=p$ for every $(x,p)\in\R^N\times\R^N$.
\medskip

\begin{lemma}\label{regular} Assume $\psi:\R^N \to \R$ to be   $\kappa_0$--Lipschitz--continuous
in $\R^N$, for some $\kappa_0 >0$ and, in addition, of class
$\CC^{1,1}$ in a neighborhood of $B_1(x_0)$, for some $x_0 \in
\R^N$. Then there exist   $t_0 >0$ and $A>0$, solely depending on
$\kappa_0$ and the Lipschitz constant of $D \psi$ in ${B_1(x_0)}$,
 such that, for every $t \in [0,t_0]$, the following properties hold:
 \begin{enumerate}
 \item[\em (i)] if $y\in B_{1/2}(x_0)$, then it is the unique  point in $\R^N$ that is optimal for $(\Tm \psi)(x)$ with  $x=\pi_1\comp\phi_t^H(y,D(\psi(y))$;\smallskip
     \item[\em (ii)]  $T^-_t \psi \in \CC^{1,1}(B_{1/2}(x_0))$\quad and \quad $\D{Lip}\big(D_x(\Tm \psi);B_{1/2}(x_0 ))\leq A$.
 \end{enumerate}
\end{lemma}

\begin{dimo}
For every $t\geq 0$, we define a map $R_t:B_1(x_0)\to\R^N$ by setting
\[
R_t(y)=\pi_1\comp\phi_t^H(y,D\psi(y))\qquad\hbox{for every $y\in B_1(x_0)$.}
\]
We claim that we can choose $t_0>0$, only depending on $\kappa_0$ and on the Lipschitz constant of $D\psi$ in $B_1(x_0)$,  such that the map $R_t - I$ is a contraction on $B_1(x_0)$ for every $t\in [0,t_0]$.
To this aim, fix $y_1$ and
$y_2$ in  $B_1(x_0)$ and denote by $(\xi_i,\eta_i)$, $i=1,2$, the
characteristic taking the value $(y_i,D\psi(y_i))$ at $0$. We
have
\[R_t(y_i) -y_i= \int_0^t \frac{\partial H}{\partial p} (\xi_i,\eta_i) \dd s\]
and consequently
\begin{equation}\label{regular01}
    |(R_t(y_1) -y_1)-(R_t(y_2) -y_2)| \leq \int_0^t
    \left |\frac{\partial H}{\partial p} (\xi_1,\eta_1)- \frac{\partial H}{\partial p}(\xi_2,\eta_2) \right | \dd s.
\end{equation}
Now we recall that $H$ is uniformly superlinear and locally bounded in $p$, uniformly in $x$, and
stays constant on characteristics. From the inequality $|D\psi(y_i)|\leq\kappa_0$ we derive that
$(\xi_i(s),\eta_i(s))$ is contained $\R^N \times B_\rho$ for some $\rho$ only depending on $\kappa_0$ and on the functions $\alpha$, $\beta$ appearing in assumption (H3). Let us denote by $\ell$ the positive constant $L_R$ given by hypothesis (T3) with $R:=\rho$. Resuming our estimate, from
\eqref{regular01} we get
\[
    |(R_t(y_1) -y_1)-(R_t(y_2) -y_2)| \leq \ell \, \int_0^t
    |(\xi_1,\eta_1)- (\xi_2,\eta_2) | \dd s.
\]
Let us denote by $\lambda$ a Lipschitz constant for $D\psi$ in $B_1(x_0)$. By Gronwall inequality
\begin{eqnarray}\label{Lipschitz flow}
  |(\xi_1(s),\eta_1(s))-(\xi_2(s),\eta_2(s)) |
  &\leq&
  e^{\ell\,s}\,|(y_1,D\psi(y_1))- (y_2,D\psi(y_2))|\\
  &\leq&
   e^{\ell\,s}\left(\sqrt{1 + \lambda^2}\right)\, |y_1-y_2|\qquad\hbox{for $s\in [0,t]$,}\nonumber
\end{eqnarray}
hence
\begin{eqnarray}\label{eq contraction}
  |(R_t(y_1) -y_1)-(R_t(y_2) -y_2)| \leq \left(\sqrt{1 + \lambda^2}\right)(e^{\ell\,t}-1)\, |y_1-y_2|.
\end{eqnarray}
Let us choose $t_0>0$ such that
\[
 (e^{\ell\,t_0}-1)\sqrt{1 + \lambda^2}< \frac 1 2,\qquad \ell\,t_0< \frac 1 4,\qquad t_0 R(\kappa_0) < \frac 14
\]
and fix $t\in [0,t_0]$. From \eqref{eq contraction} we derive that $R_t-I$ is a contraction, with Lipschitz constant less than $1/2$. By the triangular inequality we get
\[
\frac 12 |y_1-y_2|\leq |R_t(y_1)-R_t(y_2)|\leq \frac 3 2|y_1-y_2|\qquad\hbox{for every $y_1,\,y_2\in B_1(x_0)$,}
\]
showing in particular that $R_t$ is injective in $B_1(x_0)$.

Let us prove {\em (i)}. We first observe that the equality
\[
 (\Tm \psi)(x)=\psi(y)+h_t(y,x)\qquad\hbox{for $x\in B_{3/4}(x_0)$}
\]
implies $y\in B_1(x_0)$ and $R_t(y)=x$,
in view of Propositions \ref{prop B4} and \ref{prop minimizer}, in particular such a point $y$ is unique by injectivity of $R_t$.
On the other hand, by the choice of $\ell$ and the inequality $\ell t<1/4$, we easily see that the point $x=R_t(y)$ belongs to $B_{3/4}(x_0)$ if $y\in B_{1/2}(x_0)$, in particular, as previously noticed,  $y$ is the unique point in $\R^N$ that that is optimal for $(\Tm \psi)(x)$.

Let us prove {\em (ii)}. By the previous argument get that $B_{3/4}(x_0)\subseteq R_t(B_1(x_0))$, so the map
\[
 R_t^{-1}:B_{1/2}(x_0)\to B_1(x_0).
\]
is well defined and Lipschitz--continuous.
From this and \eqref{Lipschitz flow} we infer that the function
\[
x \mapsto \pi_2 \comp\phi_t^H\Big(R_t^{-1}(x),D(\psi(R_t^{-1}(x)))\Big)
\]
is $A$--Lipschitz continuous in $B_{1/2}(x_0) $, where $A$ is a constant only depending on $t_0$ and on the Lipschitz constant of $D\psi$ in $B_1(x_0)$. In addition, by Proposition \ref{prop minimizer},
its images belong, for any $x$, to $D^+(T^-_t \psi)(x)$.
This means that we have constructed an $A$--Lipschitz continuous selection of
$D^+(T^-_t \psi)(x)$ in $B_{1/2}(x_0)$. Now $T^-_t \psi$ is semiconcave, hence by exploiting basic
properties of superdifferential of semiconcave  functions we conclude that
$T^-_t \psi \in \CC^{1,1}(B_{1/2}(x_0))$ and $\D{Lip}\big(D_x(\Tm \psi);B_{1/2}(x_0 ))\leq A$.
\end{dimo}

\medskip

\begin{prop}\label{semi} Let $w$ be a critical subsolution which is $K$--semiconvex
in $B_1(x_0)$ for some $x_0 \in \R^N$. Then
 we can find $t_0 >0$, only depending  on $K$, such that
 \begin{eqnarray*}
    T^-_{t} w  \in \CC^{1,1}(B_{1/4}(x_0 ))  \quad\hbox{ for
any $t \in (0,t_0]$.}
 \end{eqnarray*}
Moreover, $\D{Lip}\big(D_x(\Tm w);B_{1/4}(x_0 ))\leq B_t$\quad for some $B_t$ only depending on $t$ and $K$.
\end{prop}

\begin{dimo}
By semiconvexity assumption on  $w$ we have that, for
any $y \in B_1(x_0)$ and $p \in D^-w(y)$, the paraboloid
\[
\psi_{y,p}(x)=w(y) + \langle p, x-y\rangle - K \, |x-y|^2,\qquad x\in B_1(x_0)
\]
is a subtangent in $B_1(x_0)$ to $w$ at $y$. We suitably extend it outside
$B_1(x_0)$ to obtain a family of $\kappa_0$--Lipschitz continuous subtangents to $w$ in $\R^N$,
i.e. satisfying $\psi\leq w$ in $\R^N$, which are in addition $C^\infty$
and possess $K$--Lipschitz differentials in $B_1(x_0)$.
Here $\kappa_0$ is a constant greater than or equal to $\kappa$. We denote
by $F$ this family of modified paraboloids.

We are then in the framework of Lemma
\ref{regular}, and consequently  there exist two common values $t_0>0$ and $A>0$ such that, for every $t\in [0,t_0]$, all
functions of $F$ satisfy items {\em (i)} and {\em (ii)} of the lemma.
Up to choose a smaller $t_0$ if
necessary, we can also assume
\begin{equation}\label{semi1}
    t_0 \,R(\kappa_0) < \frac 1 4.
\end{equation}
Let us fix $t\in [0,t_0]$. We want to show that
\begin{equation}\label{tesi semi}
(\Tm w)(x)=\sup_{\psi\in F}(\Tm\psi)(x)\qquad\hbox{for every $x\in B_{1/4}(x_0)$}.
\end{equation}
By monotonicity of $\Tm$, it is clear that $\Tm w\geq\sup_{\psi\in F}\Tm\psi$ on $\R^N$. Let us show equality on
$B_{1/4}(x_0)$. Fix $x\in B_{1/4}(x_0)$ and let $y\in\R^N$ be optimal for $(\Tm w)(x)$. Then $y\in B_{1/2}(x_0)$ in view of Proposition \ref{prop B4} and of \eqref{semi1}. By Proposition \ref{prop minimizer}, there furthermore exists
$p\in D^-w(y)$ such that
\[
x=\pi_1\comp\phi_t^H(y,p).
\]
By the very definition of $F$, there is $\psi \in F$ which is a subtangent
to $w$ at $y$ and satisfies $D\psi(y)=p$. Hence, by Lemma \ref{regular}--{\em (i)}, $y$ is also
optimal  for $(\Tm\psi)(x)$, so that
\[
(\Tm w)(x) = (\Tm \psi)(x),
\]
proving \eqref{tesi semi}.

Summing up, we have showed that $T^-_t w$ is the supremum in $B_{1/4}(x_0)$ of a family of $A$--semiconvex functions, so it is $A$--semiconvex.
Now $\Tm w$ is also $K_t$--semiconcave for some $K_t$ only depending on $t$, in view of Proposition \ref{lemma ht}. We conclude that $\Tm w$ is of class $C^{1,1}$ in $B_{1/4}(x_0)$ with $\D{Lip}\big(D_x(\Tm w);B_{1/4}(x_0 ))\leq B_t:=\max\{A,K_t\}$.
\end{dimo}

\medskip

\medskip

We go back to stationary ergodic setting  and apply  to it the
information  on the action of Lax--Oleinik semigroups  gathered in
the two previous  deterministic results as well  as other
distinguished outputs of the analysis so far performed.

\medskip

\bigskip

\noindent \textbf{Proof of Theorem  \ref{sawbis}}.  Let $w\in\S$ be a strict subsolution and let $\eps>0$ be fixed.
Choose $s>0$ such that $R(\kappa) s<\eps/2$ and set
\[
v(x,\omega)=(T^+_s w)(x,\omega),\qquad (x,\omega)\in\R^N\times\Omega.
\]
According to Proposition \ref{invariance}, $v$ is still an admissible and strict critical subsolution. Furthermore,
 by Proposition \ref{lemma ht}, there exists $K>0$ such that $v(\cdot,\omega)$ is $K$--semiconvex in $\R^N$ for every $\omega\in\Omega$. In view of Proposition \ref{semi}, there exist $t>0$ and a constant $K_t$, only depending on $t$, such that $\Tm v(\cdot,\omega)$ is a function of class $C^{1,1}$ in $\R^N$ and $\D{Lip}\big(D_x(\Tm v(\cdot,\omega));\R^N\big)\leq K_t$, for every $\omega\in\Omega$. Since $t$ can be taken arbitrarily small, we choose $t>0$ such that $R(\kappa) t<\eps/2$ and we set
\[
 w_\eps(x,\omega)=(\Tm v)(x,\omega)=(\Tm\comp T^+_s w)(x,\omega),\qquad (x,\omega)\in\R^N\times\Omega.
\]
According to Proposition \ref{invariance}, $w_\eps$ is still an admissible and strict critical subsolution. From Proposition \ref{prop B4} we get
\begin{eqnarray*}
  \|w_\eps(\cdot,\omega)-w(\cdot,\omega) \|_\infty
  &\leq&
  \|\Tm\comp T^+_s w(\cdot,\omega)-T^+_s w(\cdot,\omega) \|_\infty
  +
  \|T^+_s w(\cdot,\omega)-w(\cdot,\omega) \|_\infty\\
  &\leq&
  (t+s)R(\kappa)<\eps.
\end{eqnarray*}
Last, by Proposition \ref{prop aubryconst}, for every $\omega$ in a set of probability 1 we have
\[
w_\eps(\cdot,\omega)=\Tm\comp T^+_s w(\cdot,\omega)=T^+_s w(\cdot,\omega)=w(\cdot,\omega)\qquad\hbox{on $\A(\omega)$.}
\]
\qed
\end{subsection}

\end{section}


\begin{thebibliography}{99}

\small{

\bibitem{Ar98} \textsc{M. Arisawa,} Multiscale homogenization for
first-order Hamilton-Jacobi equations. Proceedings of the workshop
on Nonlinear P.D.E., Saitama University, 1998.


\bibitem{BCD97} \textsc{M. Bardi, I. Capuzzo Dolcetta,}
{Optimal control and viscosity solutions of
Hamilton--Jacobi--Bellman equations.} With appendices by Maurizio
Falcone and Pierpaolo Soravia.  Systems \& Control: Foundations \&
Applications. Birkh\"auser Boston, Inc., Boston, MA, 1997.

\bibitem{Ba94} \textsc{G. Barles,}
{Solutions de viscosit\`e  des \'equations de  Hamilton--Jacobi}.
Math�matiques \& Applications, { 17}.  Springer--Verlag, Paris,
1994.
%

\bibitem{Be} \textsc{P.Bernard,}  Existence of $C^{1,1}$ subsolution of the Hamilton--Jacobi equation
on compact Manifolds.  {\em Annales scientifiques de l'ENS.} {\bf
40} (2007), 445--452.

%
\bibitem{BuGiHi98} \textsc{G. Buttazzo, M. Giaquinta, S. Hildebrandt,} One--dimensional
variational problems. An introduction. Oxford Lecture Series in
Mathematics and its Applications, {\bf 15}. The Clarendon Press,
Oxford University Press, New York, 1998.
%
%
\bibitem{CaSi} \textsc{P. Cannarsa, C. Sinestrari,} Semiconcave
functions, Hamilton--Jacobi Equations, and Optimal Control.
Progress in Nonlinear Differential Equations and their
Applications, 58. Birkh\" auser Boston, Inc., Boston, MA, 2004.
%
%
%
%





\bibitem{Cl} \textsc{F.H. Clarke,} Optimization and nonsmooth analysis. Wiley, New York,



\bibitem{CoIt00} \textsc{G. Contreras, R. Iturriaga,} Global Minimizers of
Autonomous Lagrangians. 22nd Brazilian Mathematics Colloquium,
IMPA, Rio de Janeiro, 1999.


%
\bibitem{D1-06} \textsc{A. Davini,} Bolza Problems with
discontinuous Lagrangians and Lipschitz continuity of the value
function.  {\em SIAM J. Control Optim.}  {\bf 46}  (2007),  no. 5,
1897--1921.


\bibitem{DS12} \textsc{A. Davini, A. Siconolfi,} Weak KAM Theory
topics in the stationary ergodic setting, {\em Calc. Var. and
PDEs}, to appear.

\bibitem{DS11} \textsc{A. Davini, A. Siconolfi,} Metric
Techniques for convex stationary ergodic Hamiltonians, {\em Calc.
Var. and PDEs}, vol. 40 (2011), p. 391--421.

\bibitem{DS09} \textsc{A. Davini, A. Siconolfi,} Exact and
approximate correctors for stochastic Hamiltonians: the
1--dimensional case, {\em Math. Annalen,} Vol. 345, no. 4 (2009),
pp. 749--782.

\bibitem{DS06} \textsc{A. Davini, A. Siconolfi,} A generalized dynamical approach
to the large time behavior of solutions of Hamilton--Jacobi
equations. {\em SIAM J. Math. Anal.}, Vol. 38, , no. 2 (2006),
478--502.


\bibitem{DavZav}
\textsc{A. Davini, M. Zavidovique,} Weak KAM theoretic aspects for
nonregular commuting Hamiltonians. {\em Preprint} (2011). (ArXiv:
1102.2334)


\bibitem{Fathi} \textsc{A. Fathi,} Weak Kam Theorem in Lagrangian Dynamics.
Cambridge University Press, {\em to appear}.


\bibitem{FaFiRi}
\textsc{A. Fathi, A. Figalli, L. Rifford,} On the Hausdorff
dimension of the Mather quotient. {\em Comm. Pure Appl. Math.}
{\bf 62} (2009), no. 4, 445--500.


\bibitem{FaSic03} \textsc{A. Fathi, A. Siconolfi}, PDE aspects of
Aubry--Mather theory for continuous convex Hamiltonians. {\em
Calc. Var. Partial Differential Equations} {\bf 22},   no. 2
(2005) 185--228.


\bibitem{Is} \textsc{H. Ishii}, Almost periodic homogenization of Hamilton-Jacobi equations.
International Conference on Differential Equations, Vol. 1, 2
(Berlin, 1999),  600--605, World Sci. Publ., River Edge, NJ, 2000.


\bibitem{JiKoOl} \textsc{V.V. Jikov, S.M. Kozlov, O.A.
Oleinik}, Homogenization of differential operators and integral
functionals. Translated from the Russian by G.A. Yosifian.
Springer-Verlag, Berlin, 1994.



\bibitem{LPV}
\textsc{P.L. Lions, G. Papanicolau, S.R.S. Varadhan,}
Homogenization of Hamilton--Jacobi equations, unpublished preprint
(1987).

\bibitem{LiSou03} \textsc{P.L. Lions, P.E. Souganidis,}
Correctors for the homogenization of Hamilton-Jacobi equations in
the stationary ergodic setting. {\em Comm. Pure Appl. Math.} {\bf
56} (2003),  no. 10, 1501--1524.

\bibitem{Molchanov} \textsc{I. Molchanov}, Theory of random sets. Probability and its
Applications (New York). Springer-Verlag London, Ltd., London,
2005.

%
%

\bibitem{ReTa00} \textsc{F. Rezakhanlou, J. E. Tarver,} Homogenization for
stochastic Hamilton-Jacobi equations.  {\em Arch. Ration. Mech.
Anal.} {\bf 151}  (2000),  no. 4, 277--309.



\bibitem{Ru90} \textsc{W. Rudin,} Walter Fourier analysis on groups.
Reprint of the 1962 original. Wiley Classics Library. A Wiley-Interscience Publication.
John Wiley \& Sons, Inc., New York, 1990.



\bibitem{Souga99} \textsc{P. E. Souganidis},
Stochastic homogenization of Hamilton-Jacobi equations and some
applications.  {\em Asymptot. Anal.}  {\bf 20}  (1999),  no. 1,
1--11.

\bibitem{Tsi01} \textsc{B. Tsirelson}, Filtrations of random processes
in the light of classification theory. I. A topological zero-one
law. {\em Preprint} (2001) (ArXiv: math/0107121).

%
%

 }


\end{thebibliography}
\end{document}